\documentclass[12pt]{amsart}
\usepackage{amssymb}
\usepackage[all]{xy}
\calclayout 
\newtheorem{dummy}{anything}[section]
\newtheorem{theorem}[dummy]{Theorem}
\newtheorem*{thma}{Theorem A}
\newtheorem*{thmb}{Theorem B}
\newtheorem{lemma}[dummy]{Lemma}
\newtheorem{proposition}[dummy]{Proposition}
\newtheorem{corollary}[dummy]{Corollary}

\theoremstyle{definition}%%Change Theoremstyle
\newtheorem{definition}[dummy]{Definition}
 \newtheorem{example}[dummy]{Example}
 
 \newtheorem{remark}[dummy]{Remark}

\newcommand
{\eqncount}{\setcounter{equation}{\value{dummy}}%
\addtocounter{dummy}{1}}
\newcommand{\cC}{\mathcal C}

\newcommand{\cS}{\mathcal S}
\newcommand{\cD}{\mathcal D}
\newcommand{\cE}{\mathcal E}

\newcommand{\bH}{\mathbf H}
\newcommand{\bZ}{\mathbf Z}

\newcommand{\bA}{\mathbb A}
\newcommand{\bB}{\mathbb B}
\newcommand{\bC}{\mathbb C}
\newcommand{\bD}{\mathbb D}
\newcommand{\bE}{\mathbb E}
\newcommand{\bG}{\mathbb G}
\newcommand{\bR}{\mathbb R}

\newcommand{\fS}{\mathfrak E}

\newcommand{\cy}[1]{\bZ/{#1}}

\newcommand{\wK}{\widehat K}
\newcommand{\wL}{\widehat L}

\newcommand{\bd}{\partial}
\newcommand{\vv}{\, | \,}
\newcommand{\trf}{tr{\hskip -1.8truept}f}

\newcommand{\mmatrix}[4]{\left (\vcenter
{\xymatrix@C-2pc@R-2pc{#1&#2\\#3&#4} } \right )}

\DeclareMathOperator{\Mod}{\,mod}
\DeclareMathOperator{\rank}{rank}

\DeclareMathOperator{\Image}{Im}

\newcommand{\newt}{\tau^{NEW}}
\newcommand{\newtbar}{\bar{\tau}^{NEW}}

\DeclareMathOperator{\sign}{sign}
\newcommand{\mf}{}
\newcommand{\rhL}{L_{{\mf r}h}}
\newcommand{\rL}{L_{\mf r}}
\newcommand{\eL}{L_{\mf ev}}
\newcommand{\rsL}{L_{{\mf r}s}}
\newcommand{\esL}{L_{{\mf ev}\,s}}
\newcommand{\ehL}{L_{{\mf ev}\,h}}

\newcommand{\di}{}
\newcommand{\tauiso}{\tau^{iso}}
\newcommand{\spt}{\textstyle}
\newcommand{\FC}{{\mathcal C}^{fil}}
\newcommand{\SD}{\mathbb{SD}}
\newcommand{\SPD}{\mathbb{SPD}}
\newcommand{\fildual}{F^{dual}}

\begin{document}

\title[Signatures of Fibre Bundles]
{The Signature of a Fibre Bundle is Multiplicative Mod 4}
\author{I.~Hambleton}
\thanks{{\hskip -12truept}Research partially supported by
NSERC Research Grant A\,4000, the British Association Travel Fund,
and the Edinburgh Mathematical Society. The authors also wish to
thank the Max Planck Institut f\"ur Mathematik, Bonn, for its
hospitality and support.}
\address{Department of Mathematics \& Statistics
 \newline\indent
McMaster University
 \newline\indent
Hamilton, ON  L8S 4K1, Canada} \email{ian{@}math.mcmaster.ca}
%second author
\author{A.~J.~Korzeniewski}
\address{School of Mathematics
 \newline\indent
University of Edinburgh
 \newline\indent
Edinburgh EH9 3JZ, United Kingdom}
\email{A.J.Korzeniewski@sms.ed.ac.uk}
%third author
\author{A.~A.~Ranicki}
\address{School of Mathematics
 \newline\indent
University of Edinburgh
 \newline\indent
Edinburgh EH9 3JZ, United Kingdom} \email{A.Ranicki@ed.ac.uk}

\begin{abstract}\noindent
We express the signature modulo 4 of a closed, oriented,
$4k$-dimensional $PL$ manifold as a linear combination of its
Euler characteristic and the new absolute torsion invariant
defined in Korzeniewski \cite{ak1}. Let $F \to E \to B$ be a $PL$
fibre bundle, where $F$, $E$  and $B$ are closed, connected, and
compatibly oriented $PL$ manifolds. We give a formula for the
absolute torsion of the total space $E$ in terms of the absolute
torsion of the base and fibre, and then  combine these two results
to prove that the signature of $E$ is congruent modulo 4 to the
product of the signatures of $F$ and $B$.
\end{abstract}

%\date{Version 4, April 13, 2005}
\maketitle

 The \emph{signature} $\sign(M) \in \bZ$ of a closed, oriented
$n$-manifold $M^n$ is the index of its cup product form when $n
\equiv 0 \Mod 4$, and zero otherwise.  For an orientable
differentiable fibre bundle, Chern, Hirzebruch and Serre
\cite{chs1} proved that $\sign(E) = \sign(F)\cdot\sign(B)$
provided that the fundamental group $\pi_1(B)$ acts trivially on
the real cohomology $H^*(F;\bR)$ of the fibre.  In general, the
signature is \emph{not} multiplicative for differentiable fibre
bundles:  the first examples were constructed by Kodaira
\cite{kodaira1}, Atiyah \cite{atiyah1}, and Hirzebruch
\cite{hirzebruch1}.  These examples occur in the lowest possible
dimension where $\dim B = \dim F =2$.  The total space $E$ of the
bundle in these examples has dimension $4$ and non-zero signature,
but $\sign(B) = \sign(F) =0$ by convention. Surface bundles over
surfaces were studied in detail by W.  Meyer \cite{meyer1} and H.
Endo \cite{endo1}.  For such bundles, the signature of the total
space was shown to be divisible by $4$.

 The non-multiplicative behaviour of the signature in fibre bundles is
explained by the characteristic class formula for $\sign(E)$ due to Atiyah
\cite{atiyah1}, which contains a contribution from the
action of $\pi_1(B)$ on the cohomology of the fibre.  On the other
hand, W.  Neumann \cite{neumann1} showed that the signature is
multiplicative in many cases where the structural group of the bundle
is discrete, and gave a Witt group interpretation of the obstruction to
multiplicativity under this assumption.

Stephan Klaus and Peter Teichner
\cite{kt1}
have conjectured that
the signature is multiplicative modulo 8, provided that $\pi_1(B)$ acts
trivially on the mod 2 cohomology $H^*(F;\cy 2)$ of the fibre.  This
 suggests studying the conditions under which the signature is
multiplicative modulo other powers of 2.  The first case is easy since
$\sign(E)\equiv \chi(E) \Mod 2$, and the Euler characteristic $\chi(E)
=\chi(F)\cdot\chi(B)$ is multiplicative.  The following is our main
result.

\begin{thma}
Let $F \xrightarrow{q} E \xrightarrow{p} B$ be a $PL$
fibre bundle of closed, connected, compatibly oriented $PL$ manifolds.
Then
$$\sign(E) \equiv \sign(F)\cdot \sign(B) \Mod 4\ .$$
\end{thma}

We will use the new absolute torsion invariant $\newt(E) \in
K_1(\bZ[\pi_1(E)])$ defined and investigated in \cite{ak1}. The
new invariant is defined more generally, and has better properties
than the ``round torsion" invariant introduced in \cite{ra13}.  We
show in Theorem \ref{signmod4} that the value of $\sign(E)\Mod 4$
can be expressed in terms of $\newt(E)$ and $\chi(E)$.  Our main
result follows from this fact, and a general formula for
$\newt(E)$ in terms of $\newt(F)$ and $\newt(B)$.
\begin{thmb}\label{bigtheorem}
Let $F \xrightarrow{q} E \xrightarrow{p} B$ be a $PL$ fibre bundle
of closed, connected, compatibly oriented $PL$ manifolds.
  Then, if $n=\dim E$,
$$\newt(E) = p^!(\newt(B)) + \chi(B)q_{\ast}(\newt(F))
\in \widehat H^n(\cy{2};K_1(\bZ\pi_1(E)))\ .$$
\end{thmb}
In this statement the maps $p^!$ and $q_{\ast}$ are the transfer
and push-forward maps associated to the fibre bundle (see
\cite{lueck1}, \cite{lueck2}, and \cite{lueck-ranicki1} for the
surgery transfer).

In Section \ref{zero} we describe some of the results of \cite{ak1},
including the definition and main properties of the absolute torsion
invariant. In Section \ref{one} we use the new invariant to
 correct the definition of the round \emph{simple} symmetric
$L$-groups of \cite{hrt1}, so they fit in the long exact sequence
$$\dots \to L^n_{rs}(R) \xrightarrow{\hphantom{xx}}
L^n_{rh}(R) \xrightarrow{\newt} \widehat H^n(\cy{2};K_1(R))
 \xrightarrow{\hphantom{xx}} L^{n-1}_{rs}(R)\to \dots$$
 Section \ref{two} establishes the
relationship between the signature mod 4 and the new torsion
invariant of a $4k$-dimensional
Poincar\'e space (see Theorem \ref{signmod4}).
Sections \ref{polyhedra} and  \ref{four} study
 the set $\cS_\bullet(X)$ of \emph{pointed torsion structures}
(PTS) on manifolds. In Theorem \ref{manifold} we give a formula expressing the
absolute torsion invariant for a closed manifold $X$ in terms of a pairing
$$\Phi\colon \cS_\bullet(X)\times \cS_\bullet(X) \to
K_1(\bZ[\pi_1(X,x_0)])$$
After this preparation, we turn to the
geometric applications of absolute torsion for fibre bundles.In
Section \ref{transport} we define the fibre transport and transfer
map $p^!$ associated to a fibre bundle, and establish some useful
properties. Sections \ref{five} and \ref{six} deal with pointed
torsion structures
 on  fibre bundles, and the $\Phi$-pairing for the total space $E$.
Section \ref{thmB} uses the results of Section \ref{six} to prove Theorem B, and then
Section \ref{seven} uses Theorem B to prove Theorem A.
Section \ref{eight} is a brief discussion of fibrations of Poincar\'e spaces.

The last part of the paper (Section \ref{three})  contains the
algebra needed to describe the absolute torsion of signed filtered
chain complexes, such as the chain complexes of the total spaces
of fibre bundles. This may be of independent interest, and this
part can be read independently of the rest of the paper.
 We give a purely algebraic treatment of
the absolute torsion of filtered chain complexes: the main result
is the Invariance Theorem \ref{invariance1} identifying the
torsion of a contractible filtered chain complex with the torsion
of the contractible chain complex (in the derived category) of
filtration quotients.

\tableofcontents

We are indebted to Stephan Klaus and Peter Teichner for informing us
about their results, and pointing out this related question to the
first author in discussions at Oberwolfach in July 2001.

\section{ Signed $K$-theory\label{zero}}

In this section we recall the algebraic $K$-groups $K_0(R)$,
$K_1(R)$ of a ring $R$, and the algebraic theory of torsion of
Ranicki \cite{ra13}, \cite{ra16}, \cite{ra18} and Korzeniewski
\cite{ak1}. The absolute torsion invariant $\newt(f) \in K_1(R)$
of \cite{ak1} is based on the definition of ``round torsion" given
in \cite[7.20(ii)]{ra18}, but the sign correction terms are
different.

The \emph{class group} $K_0(\bA)$ of an additive category $\bA$ is
the abelian group with one generator $[M]$ for each object $M$ in
$\bA$, and relations
\begin{itemize}
\item[(i)] $[M]=[M']$ if $M$ is isomorphic to $M'$,
\item[(ii)] $[M \oplus N]=[M]+[N]$ for objects $M,N$ in $\bA$.
\end{itemize}

The \emph{isomorphism torsion group} $K^{iso}_1(\bA)$ is the abelian group with one generator
$\tau^{iso}(f)$ for each isomorphism $f\colon M \to N$ in $\bA$, and relations
\begin{itemize}
\item[(i)] $\tau^{iso}(gf)=\tau^{iso}(f)+\tau^{iso}(g)$ for isomorphisms $f\colon M \to N$, $g\colon N \to P$.
\item[(ii)] $\tau^{iso}(\begin{pmatrix} f & e \\ 0 & f'\end{pmatrix})
=\tau^{iso}(f)+\tau^{iso}(f')$
for isomorphisms
$f\colon M \to N$, $f'\colon M' \to N'$ and any morphism $e\colon M' \to N$.
\end{itemize}

The class group $K_0(R)$ and the  isomorphism torsion group $K^{iso}_1(R)$
of a ring $R$ are the class and torsion groups of the additive category
$\bA(R)$ of based f.g. (= finitely generated) free $R$-modules.
To obtain values in $K_1(R)$, we compose with the
 split surjection $K^{iso}_1(R) \to K_1(R)$
given by \cite[Prop.  5.1]{ra13}, and let $\tau^{iso}(f) \mapsto \tau(f)\in K_1(R)$.
In particular, there are defined isomorphisms
$$\begin{array}{l}
K_0(\bZ) \xymatrix{\ar[r]^-{\cong}&} \bZ~;~[M] \mapsto {\rm rank}_{\bZ}\,M~,\\[.6ex]
K_1(\bZ) \xymatrix{\ar[r]^-{\cong}&} \bZ/2~;~
[M,f] \mapsto \tau({\det}(f))~.
\end{array}$$
The \emph{reduced class group} $\widetilde{K}_0(R)$ and the
\emph{reduced torsion group} $\widetilde{K}_1(R)$ are defined by
$$\widetilde{K}_i(R)={\rm coker}(K_i(\bZ) \to K_i(R))~~(i=0,1)~.$$
We shall assume that $R$ is such that the rank of f.g.  free $R$-modules is
well-defined, so that $K_0(R)=K_0(\bZ) \oplus \widetilde{K}_0(R)$.

We shall consider only \emph{finite} chain complexes over $R$, meaning
finite dimensional, positive chain complexes of f.g.  free $R$-modules.
The \emph{Euler characteristic} of a finite $R$-module chain complex
$C$ is the chain homotopy invariant
$$\chi(C)=\sum\limits^{\infty}_{r=0}(-1)^r \rank_R C_r \in \bZ~.$$
Our assumptions on $C$ mean that $C_r = 0$ if $r <0$ or $r > N$ for
some integer $N$, so this sum is defined.
A finite chain complex $C$ is \emph{round} if $\chi(C)= 0 \in \bZ$.

The \emph{torsion} $\tau(C)\in K_1(R)$ of a contractible finite
based $R$-module chain complex
$$C: \ 0  \to \dots \to C_{r+1} \xymatrix{\ar[r]^-{d}&} C_r
\xymatrix{\ar[r]^-{d}&} C_{r-1} \to \dots \to 0$$ was
originally defined by Whitehead \cite{whitehead} to be
$$\tau(C)=\tau(d+\Gamma) \in K_1(R)~,$$
with the isomorphism
$$d+\Gamma\colon  C_{odd}=\sum\limits_{r~{\rm odd}}C_r\xrightarrow{\hphantom{xxx}}
C_{even}=\sum\limits_{r~{\rm even}}C_r$$
defined using any contraction $\Gamma\colon 1 \simeq 0\colon C \to C$.
The value of $\tau(C)$ is independent of the choice of $\Gamma$.

The torsion of a chain equivalence $f\colon C \to D$ of finite based
$R$-module chain complexes can be defined by
$$\tau(f)=\tau(\cC(f)) \in K_1(R)$$
with $\cC(f)$ the algebraic mapping cone
$$d_{\cC(f)}=\begin{pmatrix} d_D & (-)^{r-1}f \\ 0 & d_C
\end{pmatrix}\colon \cC(f)_r=D_r \oplus C_{r-1} \to
    \cC(f)_{r-1}=D_{r-1} \oplus C_{r-2}~.$$
However, as noted in
\cite{ra13}, this definition of torsion only has good sum and
composition properties modulo the sign subgroup
$${\rm im}(K_1(\bZ) \to K_1(R))=\{\tau(\pm 1)\} \subseteq K_1(R)~.$$
The new absolute torsion uses the notion of a \emph{signed} chain
complex \cite[Def.  5]{ak1}.  This is a pair $(C, \eta_C)$, where
$C$ is as above and $\eta_C \in \Image (K_1(\bZ) \to K_1(R))$. The
image of $K_1(\bZ) = \cy 2$ in $K_1(R)$ is given by $a\mapsto
a\cdot\tau(-1)$, but the symbol
$$\tau(-1): =~ \tau(-1\colon R \to R) \in K_1(R)$$
will usually be suppressed to simplify the notation.

If $(C, \eta_C)$ is a signed chain complex with $C$ contractible,
then its absolute torsion is defined as
$$\newt(C,\eta_C) = \tau(C) + \eta_C\ .$$
If $(C,\eta_C)$ is a signed complex, then the suspension
$SC$  of $C$, where $(SC)_r = C_{r-1}$, has the sign
$\eta_{SC} = -\eta_C$. The sum $(C \oplus D, \eta_{C\oplus D})$
of two signed complexes
involves some more sign terms from $K_1(\bZ) = \cy 2$.
Let
$$\epsilon(M, N) : = \rank_R(M)\cdot \rank_R(N) \Mod 2$$
for any free $R$-modules $M$ and $N$, and
 define \cite[p. 213]{ra13}:
$$\beta(C,D) = \sum_{i>j}\left (\epsilon(C_{2i},D_{2j})
- \epsilon(C_{2i+1},D_{2j+1})\right )
                        \in \Image (K_1(\bZ)\to K_1(R))\ .$$
This is just the difference of the torsions of the
permutation isomorphisms $(C\oplus D)_{even} \to C_{even}\oplus
D_{even}$ and $(C\oplus D)_{odd} \to C_{odd}\oplus
D_{odd}$. We define
$$\eta_{C\oplus D} = \eta_C +\eta_D - \beta(C,D) + \rank(C_{odd})\cdot \chi(D)\ .
$$ From now on,  we will usually denote a signed chain complex
$(C,\eta_C)$ just by $C$, even though all the formulas will
involve these signs. It turns out that the absolute torsion of a
Poincar\'e complex $C$ is independent of the choice of sign
$\eta_C$, and this justifies our abbreviated notation (see
\cite[Prop. 26.6]{ak1}).

The absolute torsion of a chain equivalence $f\colon C \to D$ of
finite based signed complexes is defined to be
$$\newt(f) = \newt(\cC(f)) \in K_1(R)~,$$
with $\eta_{\cC(f)} = \eta_{D\oplus SC}$.
The extra sign terms needed (in comparison with the invariant of \cite{ra13}) can
be seen from the formula
$$\newt(f) = \tau(C(f)) - \beta(D, SC) + \rank_R(D_{odd})\cdot
    \chi(SC) + \eta_D - \eta_C \in K_1(R)$$
given in \cite[Lemma 12]{ak1}.

For chain equivalences of signed complexes,
the absolute torsion invariant is a chain homotopy invariant, which is additive under compositions
 and direct sums (by an appropriate modification of the proof
 of \cite[Prop. 4.2, 4.4, 4.5]{ra13}).

\begin{proposition}[{\cite[Prop. 13]{ak1}}]\label{addit}
If $f\colon C \to D$ and $g\colon D \to E$  are chain equivalences
of  finite based signed chain complexes, then
$$\newt(g\circ f) = \newt(f) + \newt(g) \in K_1(R)\ .$$
 If $f\colon C \to D$ and $f'\colon C' \to D'$ are chain equivalences
of  finite based signed chain complexes, then
$$\newt(f\oplus f') =
\newt(f) +\newt(f')  \in K_1(R)\ .$$
\end{proposition}

The formula for the
absolute torsion is a modification of the definition in \cite[p.
223, 226]{ra13}, which gave an invariant with values in $K_1^{iso}(\bA)$ for chain complexes over an additive category $\bA$.  In that setting, the signs $\eta_C$ are in the
image of the skew-symmetric pairing
$$\epsilon\colon K_0(\bA) \otimes K_0(\bA) \to K^{iso}_1(\bA)$$
defined in \cite[Prop.  2.2]{ra13}:
for any two objects $M$, $N$ of $\bA$,
$$\epsilon([M], [N]) = \tauiso(M\oplus N \to N \oplus M)\in  K^{iso}_1(\bA)$$
is the torsion of the interchange map. See \cite[Lemma 7]{ak1}
for the properties of this pairing.

This extra generality will be useful in Section
\ref{three} where we discuss the absolute torsion of filtered
chain complexes. A finite chain complex over $\bA$ is an object in the derived category $\bD(\bA)$, whose morphisms are the chain homotopy classes of chain maps.
\begin{definition}\label{signedderived}
The \emph{signed derived category} $\SD(\bA)$ of an additive category
$\bA$ is the additive category with objects signed complexes
$(C,\eta_C)$ in $\bA$, and morphisms the chain homotopy classes of chain
maps. An isomorphism in $\SD(\bA)$ is a chain homotopy class of chain
equivalences.
\end{definition}
If $f\colon C \to D$ is a morphism in $\SD(\bA)$, then the mapping cone $\cC(f)$ is an object in $\SD(\bA)$ under our sign conventions.
 The absolute torsion gives an  invariant in
$\newt(f)\in K^{iso}_1(\SD(\bA))$ for an isomorphism $f\colon C \to D$  in $\SD(\bA)$. We
refer to \cite{ak1} for the definition and properties of the absolute torsion in this setting, but  we note that the map
$$i_*\colon K^{iso}_1(\SD(\bA)) \to K^{iso}_1(\bA) $$
defined by $i_*\tau^{iso}(f) = \newt(f)$
is a split surjection.
 The map $i_*$ has the naturality property
$$i_*\epsilon(C,D)=
\epsilon(\chi(C),\chi(D))\in K^{iso}_1(\bA)$$
for any objects $(C,\eta_C)$, $(D,\eta_D)$ in $\SD(\bA)$. We will also
use the relation
$$i_*\epsilon(C,SD)= - i_*\epsilon(C,D)
\in K^{iso}_1(\bA)$$
and the formula
$$\eta_{C(f)} = \eta_{D\oplus SC} = \eta_D - \eta_{C} -\beta(D, SC) +
\epsilon(D_{odd}, \chi(SC))\in K_1^{iso}(\SD(\bA))$$
for the mapping cone sign in our sign calculations (see Section \ref{three}).

\section{ Round $L$-theory\label{one}}

In this section we recall the round $L$-theory of Hambleton, Ranicki
and Taylor \cite{hrt1}.  Note that the
round torsion is \emph{not} a round cobordism invariant, contrary to the assertion in \cite[p.~190]{ra18},  but the new absolute torsion of \cite{ak1} does have
this important property.  Moreover, all the results claimed in
\cite[7.22(ii)]{ra18} and \cite[p.~135]{hrt1} do hold after replacing
the round torsion with the new absolute torsion invariant.

A finite, oriented Poincar\'e duality space $X$ of dimension $n$, with
universal covering $\widetilde X$, has an associated symmetric chain
complex $(C,\varphi)$, as defined in \cite{ra10}.  Here
$C: =C(\widetilde X)$ is a $n$-dimensional finite chain complex over $R
= \bZ[\pi_1(X)]$ and $\varphi_0\colon C^{n-*} \to C_*$ is a chain
equivalence.  The definition of this symmetric structure $\varphi$ uses
the standard involution $g\mapsto g^{-1}$ on the group ring.  More
generally, for any ring $R$ with involution $\alpha$, and unit
$\varepsilon=\pm 1$, Ranicki \cite{ra10} defines
$\varepsilon$-symmetric structures on algebraic chain complexes and the
notion of a symmetric (algebraic) Poincar\'e complex.  The cobordism
group of $n$-dimensional symmetric Poincar\'e is denoted $L^n(R,
\varepsilon)$ \cite[\S 3]{ra10}.  The additive inverse is defined to be
$-(C,\varphi) : = (C, -\varphi)$.  In the rest of the paper we will
assume that $\varepsilon = +1$, and denote these groups just by
$L^n(R)$.

If we restrict to round complexes in both objects and bordisms
and $n>0$, we get the round symmetric $L$-groups $\rL^n(R)$.  The groups
$\rL^0(R)$ is the Witt group of formal differences $(M,\psi) -
(M',\psi')$ of non-singular symmetric forms over $(R,\alpha)$ with
$\rank_R M = \rank_R M'$ \cite[Prop.  2.1]{hrt1}.  We can compare the
round $L$-groups to the ordinary ones by a long exact sequence.  The
anti-automorphism on $R$ induces a $\cy 2$-action on $K_0(R)$ which is
the identity on the subgroup $K_0(\bZ)= \bZ$.

\begin{proposition}[{\cite[Prop. 3.2]{hrt1}}] For any ring  with
involution
$(R, \alpha)$,
 there is a long exact sequence
$$ \dots\to  \widehat H^{n+1}(\cy 2; K_0(\bZ))
 \to \rL^n(R) \to L^n(R)
\to \widehat H^n(\cy 2;  K_0(\bZ))\to  \dots$$
where the map
$L^n(R) \to \widehat H^n(\cy 2;  K_0(\bZ))$ is defined by
$(C,\varphi) \mapsto \chi(C)$.
\end{proposition}

We will also need to compare these $L$-groups to the cobordism groups $\eL^n(R)$
of $n$-dimensional
symmetric Poincar\'e complexes $(C, \varphi)$
with \emph{even} Euler characteristic $\chi(C) \equiv 0 \Mod 2$
for objects and bordisms:
\vskip 4truept
$$\vcenter{
\xymatrix@!C@C-3pc@R-1pc{ \widehat H^{n+1}(\cy 2; \cy 2) \ar[dr]
\ar@/^2pc/[rr]    && {\widehat H^n(\cy 2; \bZ)} \ar[dr]
\ar@/^2pc/[rr]               &&
\rL^{n-1}(R)                                                        \\
&  \eL^n(R)                \ar[dr] \ar[ur]                      &&
\widehat H^n(\cy 2; \bZ)                    \ar[dr] \ar[ur]                      \\
\rL^n(R)                  \ar[ur] \ar@/_2pc/[rr]_{}            &&
L^n(R)                   \ar[ur] \ar@/_2pc/[rr]               &&
\widehat H^n(\cy 2; \cy 2)
}}
$$
\medskip
\noindent
The map $\eL^n(R) \to \widehat H^n(\cy 2; \bZ)$ is given by
$(C,\varphi)\mapsto
\chi(C)/2$.

In order to define the \emph{torsion} of a symmetric Poincar\'e
chain complex over a ring with involution $(R, \alpha)$, we
need to work with \emph{signed} chain complexes, in which each
free chain module $C_r$ has a preferred basis, along with a single
choice of sign $\eta_C$. We denote by $\rhL^n(R)$ the
corresponding cobordism group of round, $n$-dimensional, finite,
signed symmetric Poincar\'e  chain complexes of finitely generated
based free $R$-modules. By forgetting the preferred bases and
signs we can identify $\rhL^n(R) =\rL^n(R)$. Similarly, we can
identify the based and unbased versions of the even $L$-groups
$\ehL^n(R) = \eL^n(R)$.

Now let $(C,\varphi)$ denote an $n$-dimensional signed symmetric
Poincar\'e chain complex, of finitely generated based free
$R$-modules, over a ring with involution $(R, \alpha)$.
 The absolute torsion
$$\newt(C,\varphi) : =\newt(\varphi_0\colon C^{n-*} \to C_*) \in K_1(R)$$
is now defined, as the absolute torsion of the  chain equivalence
 $\varphi_0$ with respect to the given bases on $C_*$
and the dual bases on $C^{n-*}$.
 Recall that the differential for $C^{n-*}$ has the
 sign convention
 $$d_{C^{n-*}} = (-1)^r d^*_C \colon C^{n-r} \to C^{n-r+1}\ .$$
 The dual signed chain complex $(C^{n-*}, \eta_{C^{n-*}})$ is given
 the sign
 \eqncount
 \begin{equation}
\label{eqn:dualsigned}
\eta_{C^{n-*}} = \eta_C + \beta(C,C) + \alpha_n(C) \in \Image(K_1(\bZ) \to K_1(R))
\end{equation}
 where
$$\alpha_n(C) = \sum_{r\equiv n+2,n+3 (\Mod 4)} \rank_R (C^r)
 \in \Image(K_1(\bZ)\to K_1(R))\ .$$
 With this convention, the absolute torsion $\newt(C,\varphi)$
 is independent of the initial choice of sign $\eta_C$
(see \cite[Prop. 26]{ak1}).
We say that $(C,\varphi)$ is \emph{simple}
 if $\newt(C, \varphi) = 0$.
 \begin{lemma}[{\cite[Lemma 10]{ak1}}]\label{isomorphism}
 If $\varphi_0\colon C^{n-*} \to C$ is an isomorphism, then
$$\newt(C,\varphi) = \sum_{r=0}^{n} (-1)^r\tau(\varphi_0\colon
 C^{n-r} \to C_r) + \beta(C,C) + \alpha_n(C) \in K_1(R)\ .
$$
 \end{lemma}
 \begin{proof} The absolute torsion of a chain isomorphism
$f\colon C \to D$ of signed chain complexes is just
$$\newt(f) = \sum_{r=0}^n(-1)^{r} \tau(f_r\colon C_r \to D_r) - \eta_C + \eta_D$$
according to \cite[Def.~6]{ak1}. The sign terms for $\varphi_0$
appear above.
  \end{proof}

   The absolute torsion has the
symmetry property \cite[Prop. 26.2]{ak1} (compare \cite[7.20(ii)]{ra18})
$$\newt(C,\varphi)^*= (-1)^n\newt(C,\varphi) +\frac{n(n+1)}{2}\chi(C)$$
where $\ast\colon K_1(R) \to K_1(R)$ denotes the
involution on $K$-theory
induced  by  ``$\alpha$-conjugate-transpose" of
matrices. The sign term vanishes over rings $R$ for which
   skew-hermitian forms necessarily have even rank (such as
  an integral group ring $R = \bZ[\pi]$).
Also
 $$\newt(C,-\varphi) = \newt(C,\varphi) + \chi(C)$$
 (see \cite[Prop. 26.5]{ak1}).
 Notice that the extra sign terms in both formulas vanish for round or
 even symmetric complexes.
 Moreover, the absolute
 torsion is additive under direct sums of symmetric $n$-complexes.
It follows that $\newt(C,\varphi)$ defines an additive map into
$\widehat H^n(\cy{2};K_1(R))$ for round or even symmetric complexes.

 For our applications, we will need the cobordism groups of
 simple round or even Poincar\'e complexes, and their relation to
 the $L$-groups already mentioned.
However, to give a well-defined homomorphism
$$ \eL^n(R) \to \widehat H^n(\cy{2};K_1(R))$$
we need the invariant associated to  a null-bordant symmetric
complex to be trivial.

\begin{example}\label{ex-forms}
Let $R = \bZ$ with trivial involution. A zero-dimensional
 symmetric Poincar\'e complex over $\bZ$ is just unimodular,
 symmetric bilinear form $(L, h)$ on a finitely-generated
 free abelian group. In this case,
  $$\newt(L,h) = \tau(\det h) \in \cy 2$$
   by Lemma \ref{isomorphism}.
  To compare this invariant to the determinant, we map
 $\cy 2 \to \cy 4$ by $a \mapsto 2a \Mod 4$. Then
 $2\tau(\det h) \equiv \det h -1\Mod 4$ and
 $$2\newt(L,h) \equiv \det h  -1 \Mod 4\ .$$
 Note that for any form $(L,h)$ of even rank, the absolute torsion is just
 $\tau(\det h)$. In particular, the absolute torsion of the hyperbolic
 plane $\bH$ is non-zero, so $\bH$ represents a non-trivial
 element of  $\eL^0(\bZ)$.
 \end{example}
\begin{lemma}[{\cite[Prop. 26.4]{ak1}}]\label{boundary}
If $(C,\varphi)$ is homotopy equivalent to the boundary of
 an $(n+1)$-dimensional symmetric signed complex $(D, \Phi)$, then
$$\newt(C,\varphi) =
  (-1)^{n+1}\newt(C \to \bd D)^* -\newt(C \to \bd D) +
  \frac{(n+1)(n+2)}{2}\chi(D)\ .$$
\end{lemma}
Once again, the extra sign term vanishes for round or even complexes.
\begin{corollary}\label{homotopic}
The absolute torsion $\newt(C,\varphi)$
 induces a well-defined
homomorphism $ \eL^n(R) \to \widehat H^n(\cy{2};K_1(R))$.
  \end{corollary}
The cobordism group of simple round (or even)
symmetric Poincar\'e complexes is denoted
$\rsL^n(R)$ or $\esL^n(R)$, following \cite[p. 135]{hrt1}.

\begin{proposition}\label{braid}
There is a commutative braid of exact sequences
$$\vcenter{
\xymatrix@!C@C-4pc@R-1pc{
\widehat H^{n+1}(\cy 2; \bZ)
  \ar[dr] \ar@/^2pc/[rr]^{}    &&
 \rhL^n(R)               \ar[dr] \ar@/^2pc/[rr]               &&
 \widehat H^{n}(\cy 2;  K_1(R))                                                   \\
&  \rsL^n(R)              \ar[dr] \ar[ur]                      &&
 \ehL^n(R)                \ar[dr] \ar[ur]                      \\
\widehat H^{n+1}(\cy 2; K_1(R))
 \ar[ur] \ar@/_2pc/[rr]_{}            &&
 \esL^n(R)                 \ar[ur] \ar@/_2pc/[rr]               &&
\widehat H^n(\cy 2; \bZ)
}}
$$
\medskip
\noindent
where the maps $\rhL^n(R) \to \widehat H^{n}(\cy 2;  K_1(R))$
and $\ehL^n(R) \to \widehat H^{n}(\cy 2;  K_1(R))$
are induced by the absolute  torsion $(C,\varphi) \mapsto
\newt(\varphi_0) \in K_1(R)$.
  The maps to $\widehat H^n(\cy 2; \bZ)$ are given by
  $(C, \varphi)\mapsto \chi(C)/2$.
\end{proposition}
    \begin{lemma}\label{chi4}
 Let $R$ be a ring with involution, and
$(C,\varphi)$ be an even based symmetric Poincar\'e complex
over  $R$,
with $\chi(C) \equiv 0 \Mod 4$ and $\{\newt(C,\varphi)\} =0 \in
\widehat H^n(\cy 2; K_1(R))$. Then $(C,\varphi)$
is even Poincar\'e bordant to a round simple based symmetric
Poincar\'e complex $(C',\varphi')$.
 \end{lemma}

\begin{proof}
 This follows from the comparison braid given in
Proposition \ref{braid}  relating $\rhL^n(R)$ and $\ehL^n(R)$.
\end{proof}
\begin{example} For the special case $R=\bZ$ and $n=0$, we
can substitute the calculation $\rhL^0(\bZ) = 2\bZ$ from \cite[4.2]{hrt1}:
$$\vcenter{
\xymatrix@!C@R-1pc{
0
  \ar[dr] \ar@/^2pc/[rr]^{}    &&
2\bZ              \ar[dr] \ar@/^2pc/[rr]               &&
 \cy 2                                                  \\
& 4\bZ              \ar[dr] \ar[ur]                      &&
 \ehL^0(\bZ)                \ar[dr] \ar[ur]                      \\
0
 \ar[ur] \ar@/_2pc/[rr]_{}            &&
 \esL^0(\bZ)                 \ar[ur] \ar@/_2pc/[rr]               &&
\cy 2
}}
$$
We obtain $\rsL^0(\bZ) = 4\bZ$.
Since the hyperbolic plane $\bH$ represents a non-zero element of
$\ehL^0(\bZ)$ with $\chi \equiv 2 \mod{4}$ and $\det(\bH) =-1$, we have
$\ehL^0(\bZ) = 2\bZ \oplus \cy{2}$, generated by $\langle 1\rangle
 \perp \langle 1\rangle$ and $\bH$,    and
$\esL^0(\bZ) = 2\bZ$ generated by $\langle 1\rangle\perp \langle 1\rangle$.
 A diagram chase gives $\ehL^2(\bZ)=\cy 2$ and the fact that
 $\{\newt(C,\varphi)\} = 0\in \widehat H^0(\cy 2;K_1(\bZ))$ if $\dim C \equiv 2 \Mod 4$.
\end{example}
 We will need the following calculation in a later section.
 \begin{lemma}\label{ordertwo} Let $C$ be a symmetric Poincar\'e complex over $R = \bZ[\cy 2]$.
  The absolute torsion
 $\{\newt(C,\varphi)\} = 0 \in \widehat H^0(\cy 2; K_1(\bZ[\cy 2]))$ if $\dim C \equiv
 2\Mod 4$.
 \end{lemma}
 \begin{proof} Since $\dim C \equiv 2 \Mod 4$, $C$ is an even symmetric
Poincar\'e complex, so $\newt(C,\varphi)$ defines an element of
$\widehat H^0(\cy 2; K_1(R))$.  Let $\cy 2 = \langle T \rangle$ denote
a generator of the group of order two, and we have an inclusion
$\bZ[\cy 2] \subset \bZ \oplus \bZ$ of rings given by $T \mapsto \pm
1$.  The induced map $$\widehat H^0(\cy 2; K_1(\bZ[\cy 2])) \to
\widehat H^0(\cy 2; K_1(\bZ\oplus \bZ))$$ is injective, and the
absolute torsion map $\ehL^{4k+2}(\bZ[\cy 2]) \to \widehat H^0(\cy 2;
K_1(\bZ[\cy 2]))$ composed with this injection factors through
$\ehL^2(\bZ) \to \widehat H^0(\cy 2; K_1(\bZ))$, which we have seen is the zero map.
 \end{proof}

\section{ Absolute torsion and signatures\label{two}}

Let $X$ be a finite, connected $CW$-complex, and fix a base point
$x_{0}\in X$.  We choose orientations and lifts to the universal
covering $\widetilde X$ for each cell in $X$, together with a sign
$\eta_{X}$.  The cellular chain complex $(C(\widetilde X), \eta_{X})$
is now a finitely-generated, signed, based chain complex over
$R=\bZ[\pi_1(X,x_{0})]$.  If $X$ is an oriented, finite, geometric
Poincar\'e $n$-complex, then we define $$\newt(X) = \newt(C(\widetilde
X), \varphi_{0}) \in \widehat H^{n}(\cy{2};K_{1}(\bZ[\pi_1(X,x_{0})]))$$ where
$\varphi_{0}$ is the duality map from the symmetric structure on
$C(\widetilde X)$ defined by Ranicki \cite[p.  92]{ra10}.  The
symmetric structure is a homotopy invariant of $X$.

 \begin{lemma}
 Let $X$ be an oriented, finite, Poincar\'e complex of dimension $n$.
The absolute torsion $\newt(X) \in
 \widehat H^n(\cy 2; K_1( \bZ[\pi_1(X)]))$ depends only on the
homotopy type of $X$. In particular, the
absolute torsion is independent of the choice of preferred base for
$C(\widetilde X)$.
 \end{lemma}
\begin{proof}
The absolute torsion $\newt(X)$ is independent of all the choices made
by \cite[Prop.  26]{ak1}.  \end{proof}

The \emph{symmetric signature} of $X$ is the element
$\sigma^*(X) \in L^n(\bZ[\pi_1(X,x_{0})])$ given by the
bordism class of the symmetric structure
$(C(\widetilde X), \varphi)$.
If $\chi(X) \equiv 0 \Mod 2$, then
we get an even symmetric signature
$\sigma^*_{\mf ev}(X) \in \eL^n(\bZ[\pi_1(X)])$.
In this case, $\newt(X)$ is a bordism invariant:  it is the image of
 $\sigma^*_{\mf ev}(X)$
under the homomorphism $ \eL^n(R) \to \widehat H^n(\cy{2};K_1(R))$
from Corollary \ref{homotopic}. Similarly, if
$\chi(X) =0$ the image of  $\sigma^*_{\mf r}(X) \in \rL^n(\bZ[\pi_1(X)])$
is again the absolute torsion $\newt(X)$.

\begin{example} The  manifold $ X =\bC
P^2\,\#\,\bC P^2 \,\#\, 2(S^1\times S^3)$ has a round
symmetric signature, whose absolute
torsion invariant in $\widehat H^0(\cy 2, K_1(\bZ))$
is non-zero. Its symmetric signature represents a generator
of $\rhL^4(\bZ) = 2\bZ$.
\end{example}
    By composing with the augmentation map
$\epsilon\colon \bZ[\pi_1(X)] \to \bZ$, setting $\epsilon(g) = 1$
  for all $g\in \pi_1(X,x_0)$, we get a symmetric signature
$\sigma^*(X) \in L^n(\bZ)$ for any
oriented, finite Poincar\'e complex of dimension $n$ (e.g.
any closed, oriented $n$-manifold).
  We define the \emph{reduced} absolute torsion, as the image
  $$\newtbar(X): =\epsilon_\ast(\newt(X)) \in K_1(\bZ)$$
   of the absolute torsion.
   This invariant is computed from the chain
  complex of $X$, instead of $\widetilde X$.
  Our main result depends on a mod 4 relationship between the
  ordinary signature of a manifold and the  absolute torsion of its symmetric signature over $\bZ$. This is an algebraic fact.

  \begin{theorem}\label{signmod4}
  Let $(C,\varphi)$ be a finite, based, $4k$-dimensional algebraic Poincar\'e complex over $\bZ$. Then
  $$\sign(C) \equiv 2\newt(C,\varphi) +  (2k+1)\chi(C) \Mod 4\ .$$
   \end{theorem}
   \begin{proof}
   The first step is to show that the right-hand side is an algebraic
cobordism invariant for Poincar\'e complexes over $\bZ$.  Suppose that
$(C \oplus C' \to D, \delta \varphi, \varphi\oplus\varphi')$ is an
algebraic cobordism.  Since $\dim C = 4k$, we have
    $$\newt(C,\varphi) + \newt(C',-\varphi') = (4k+1)(2k+1)\chi(D) \in K_1(\bZ)$$
  by Lemma \ref{boundary}.  However,
  $\chi(D) = \frac{1}{2}(\chi(C) + \chi(C')) \in \bZ$ and
  $\newt(C', -\varphi') = \newt(C',\varphi') + \chi(C')$, so that
  $$2\newt (C,\varphi) + 2 \newt (C',\varphi') \equiv
  (2k+1) \chi(C)  + (2k-1) \chi(C')
  \Mod 4\ .$$
  Therefore
  $$2\newt(C,\varphi)+ (2k+1)\chi(C)\equiv  2\newt(C',\varphi') + (2k+1)\chi(C')
  \Mod 4$$
   as required (since $2\chi(C) \equiv 2\chi(C') \Mod 4$).
  Next, we use the fact that every $4k$-dimensional symmetric Poincar\'e
  complex over $\bZ$ is algebraically cobordant to
  a complex $(C, \varphi)$ which is concentrated in dimension $2k$
  (see \cite[Prop.~4.5]{ra10}). In
  that case $\varphi_0$ is an isomorphism, so
  $$\newt(C,\varphi) = \tau(\det\varphi_0) +\beta (C,C) +
  \alpha_{4k}(C)$$
   by Lemma \ref{isomorphism}. But the $\beta$-term is zero
   in this case, and $\alpha_{4k}(C)$ contributes $\chi(C)$ when
   $k$ is odd and zero otherwise. We may express this as
   $\alpha_{4k}(C) = k\chi(C)\in K_1(\bZ)$.
   Therefore
   $$\newt(C,\varphi) = \tau(\det\varphi_0) + k\chi(C) \in K_1(\bZ)$$
   and the right-hand side becomes
   $$ 2\newt(C,\varphi)+ (2k+1)\chi(C) \equiv 2 \tau(\det\varphi_0)
   +\chi(C)\equiv \det\varphi_0 + \chi(C) -1 \Mod 4\ .$$
   But there is a  classical formula
$$\sign(\phi) \equiv \rank(\phi) + \det(\phi) -1 \Mod 4$$
   relating the
signature and the determinant mod 4
(see \cite[Theorem 3.5]{hnk1}),
for any  unimodular symmetric bilinear form $\phi$ over the integers.
 We can apply this to $\varphi_0$, and note that $\chi(C) =
 \rank \varphi_0$, since $C$ is concentrated in dimension $2k$.

 \end{proof}

\begin{corollary}\label{sign4}
 If $(C,\varphi)$ is a simple, round Poincar\'e complex over $\bZ$,
the signature $\sign(C) \equiv 0 \Mod 4$.
\end{corollary}

We now consider the case where $X$ is an even-dimensional manifold.
Let $\pi^2$ denote the subset of $\pi^{ab}$ given by the images of
  squares from  $\pi=\pi_1(X,x_0)$.
 \begin{proposition} \label{manifold-torsion}
 Let $X$ be a closed, oriented $PL$ manifold of
 even dimension. Then $\newt(X) \in \Image\left (\widehat H^0(\cy{2};K_1(\bZ))\oplus
 \widehat H^0(\cy{2};\pi^2) \to \widehat H^0(\cy{2};K_1(\bZ[\pi]))\right )$.
  \end{proposition}
  \begin{proof}
  Note that $\widehat H^0(\cy{2};\pi^2) = \pi^2\cap\{ \bar g \in \pi^{ab}\vv \bar g^2 =1\}$.
  Suppose that $\newt(X) = \tau(\pm g)$ for some $g\in \pi$ with
  image $\bar g\neq 1 \in \pi^{ab}$. Then the duality property
  $\newt(X) = \newt(X)^\ast$ implies that $\bar g^2=1$.
  Suppose that $\bar g \notin \pi^2$. Then there is a projection
  $j\colon \pi \to  \cy 2$ such that   $j(g) = T$
  is the generator of the quotient group $\cy 2$.
  Let $p\colon X' \to X$ denote the $2$-fold covering of $X$
  induced by $j$.
  If $\dim X \equiv 2 \Mod 4$, then $$j_*\{\newt(X)\} =
 \{ \tau(\pm T)\} \in \widehat H^0(\cy 2; K_1(\bZ[\cy 2]))$$
 which is non-zero, contrary to Lemma \ref{ordertwo}.

 If $\dim X=4k$, then the transfer
  $p^!(\newt(X))= \tau(-1)$, since for finite coverings
  the fibre bundle transfer is the same as the classical
  transfer induced by restriction (see \cite[p.~108]{lueck1}).
  However, by Theorem B it follows that $\newt(X') = p^!(\newt(X))$,
 so we have
  $\newtbar(X') \neq 0 \in K_1(\bZ)$. It
  follows that $\sign(X') \equiv 2 + 2\chi(X) \Mod 4$,
  from the relation $\chi(X') = 2\chi(X)$ and Theorem
  \ref{signmod4}. The corresponding formula for $X$
  gives $$\sign(X) \equiv 2\newtbar(X) + (2k+1)\chi(X)\Mod 4$$
  and by the Hirzebruch Signature Theorem we have
  $$\sign(X') = 2\sign(X) \equiv 2\chi(X) \Mod 4$$
  This is a contradiction, so the covering $X' \to X$ does not exist
  and $\bar g \in \pi^2$.
  \end{proof}

\section{ Absolute torsion structures on polyhedra\label{polyhedra}}

In this section we will define the notion of a (pointed) torsion structure on a polyhedron. First we recall some standard definitions
(following \cite[\S 1]{casson1})
A \emph{polyhedron} is a topological space equipped with a maximal family of $PL$ related locally finite triangulations. A \emph{cell complex} $K$ is a collection of cells $PL$ embedded in a polyhedron $P$ such that
\begin{enumerate}
\item $K$ is a locally finite covering of $P$,
\item if $\beta, \gamma \in K$ then $\bd \beta$ and $\beta\cap\gamma$
are unions of cells of $K$,
\item if $\beta, \gamma$ are distinct cells of $K$, then $\text{Int}\,
\beta\cap \text{Int}\,\gamma=\emptyset$.
\end{enumerate}
Let $|K|$ denote the underlying polyhedron of a cell complex, and use $\beta$ to denote a cell in $K$ or the subcomplex it determines. A cell complex $K'$ is a \emph{subdivision} of $K$ if
$|K'| = |K|$ and every cell of $K$ is a union of cells of $K'$.
Any two cell complexes $K'$, $K''$ with $|K'| = |K''|$ have a common subdivision
$K_0$. A base-point for $K$ is a preferred vertex. A pointed map is one which preserves given base-points.  A \emph{cellular map} $h\colon K \to L$ between cell complexes is a $PL$ map $h\colon |K| \to |L|$ such that for any cell $\beta \in K$, the image $h(\beta)$ is contained in a cell of $L$. Given any continuous map $f\colon |K| \to |L|$, there is a subdivision $h\colon K_0 \to K$ and a cellular map
$k\colon K_0 \to L$ such that $f\circ h \simeq k$. In other words, any map  is homotopic to a cellular map after subdividing the domain. A \emph{cellular homeomorphism} is a cellular map
$h\colon K \to L$ such that $h\colon |K| \to |L|$ is a 
homeomorphism. From now on, we will consider only \emph{finite} cell complexes (those consisting of only finitely many cells), or coverings of finite cell complexes.

If $K$ is a cell complex, we will say that a covering space $\widetilde
P \to P$ over $P = |K|$ is \emph{subordinate} to $K$ if every cell of
$K$ is contained in a evenly-covered neighbourhood of $P$.  In that
case, we let $\widetilde K$ be the cell complex on $\widetilde P$
induced by the covering $p\colon \widetilde P \to P$.  If $k_0\in K$ is
a base point and $\pi$ denotes the structural group of the covering
(i.e.  there is an identification $\pi = p^{-1}(k_0)$) , then the
cellular chain complex $C(\widetilde K)$ is a  chain complex of  free $R$-modules, where $R = \bZ[\pi]$. We need some additional data to get a chain complex of based modules.

Let $(K, p)$ denote a pointed finite cell complex $K$,  with $|K|=P$,
and a covering space $p\colon \widetilde P \to P$ subordinate to $K$, with structural group $\pi$.
A \emph{geometric basis} for  $(K, p)$ consists of the following data:
\begin{enumerate}
\item an ordering for the cells of $K$ compatible with the boundary partial ordering.
\item an orientation on each cell $\beta\in K$,
\item a preferred lift for each cell $\beta \in K$ to  $\widetilde K$,
\item a sign $\eta_K \in \Image (K_1(\bZ) \to K_1(\bZ [\pi]))$.
\end{enumerate}
We call $(K,p)$ a \emph{based cell complex} if it is equipped with
a geometric basis. Note that the chain complex $C(\widetilde K)$ of a based cell complex  is a signed, based, finite chain complex of finitely-generated free $\bZ[\pi]$-modules.

Let $(X,x_0)$ be a pointed space, and let $p\colon \widetilde X\to X$ be a covering space with structural group $\pi$. We fix a base-point
$\tilde x_0 \in \widetilde X$, corresponding to the identity element of $\pi$, with $p(\tilde x_0) = x_0$.

A \emph{pointed torsion structure} on $(X,p)$  is a pair $(K,c,k_0)$ such that
\begin{enumerate}
\item $K$ is a pointed finite cell complex, with base-point $k_0 \in K$,
\item $c\colon |K| \to X$ is a pointed homeomorphism, and
\item $(K, c^*(p))$ is a based cell complex, with the lift of $k_0$ mapping to $\tilde x_0$.
\end{enumerate}
The notation
$c^*(p)\colon |\widetilde K| \to |K|$ means the pull-back covering via $c$.
Usually we will suppress mentioning the base-points and use the notation $(K,c)$.
\begin{definition}
Two pointed torsion structures $(K_1,c_1)$ and $(K_2,c_2)$  on $(X,p)$ are \emph{related} by $(K_0,h,k)$ if $K_0$ is a pointed cell complex and there exists a homotopy commutative diagram
$$\xymatrix{K_1\ar[rr]^f&&K_2\cr &K_0\ar[ul]^h\ar[ur]_k&}$$
where  $h$ is a pointed cellular homeomorphism, $k$ is a pointed cellular map,  $f= c_2^{-1}\circ c_1$, and  $f\circ h \simeq k$.
\end{definition}
Unless explicitly mentioned,
 we assume that all homotopies of pointed maps are base-point preserving. Since $f=c_2^{-1}\circ c_1$ is homotopic to a cellular map
 on some subdivision $K_0$ of $K$, any two pointed torsion
 structures are related as above.
If $(K_1,c_1)$ and $(K_2,c_2)$  are related by $(K_0,h,k)$, we let
 $(K_0,p_0)$ denote the
pull-back covering $(c_1\circ h)^*(p) \cong (c_2\circ k)^*(p)$.
 \begin{lemma}
Suppose that $(K_1,c_1)$ and $(K_2,c_2)$
are related by $(K_0,h,k)$, and that
$(K_0, p_0)$ is a based cell complex.
Then the quantity
$$\newt(\xymatrix{K_1 \ar@{~>}[r] & K_2}) :=
\newt(C(\widetilde K_0) \xrightarrow{k_*} C(\widetilde K_2))
-\newt(C(\widetilde K_0) \xrightarrow{h_*} C(\widetilde K_1))$$
 in $K_1(\bZ[\pi])$ is independent of the choice of $(K_0, h, k)$,  and independent of the choice of geometric basis for $(K_0,p_0)$.
\end{lemma}
\begin{proof}
If we have two choices $(K_0, h, k)$ and $(K'_0, h', k')$, there is a homotopy commutative diagram
$$\xymatrix{&K'_0\ar[dl]_{h'}\ar[dr]^{k'}&\cr
K_1\ar[rr]^f&&K_2\cr &K_0\ar[ul]^h\ar[ur]_k}$$
and we apply the chain homotopy invariance of the absolute torsion, and the composition formula of  Proposition \ref{addit}, to the composite
$C(h)^{-1}\circ C(h')\simeq C(k)^{-1}\circ C(k')$ of chain equivalences
(and chain homotopy inverses).
\end{proof}
\begin{definition}
Two pointed torsion structures $(K_1,c_1)$, $(K_2,c_2)$ on $(X,p)$
are \emph{equivalent} if they are related by a based cell complex
$(K_0,p_0)$ and
$$\newt(\xymatrix{K_1 \ar@{~>}[r] & K_2})=0 \in K_1(\bZ[\pi])\ .$$
\end{definition}
\begin{lemma} This definition gives an equivalence relation on the
set of pointed torsion structures on $(X,p)$.
\end{lemma}
\begin{proof} The given relation is reflexive and symmetric. If
$(K_1,c_1)$, $(K_2,c_2)$ are related by $(K_0,h,k)$,
and $(K_2,c_2)$, $(K_3,c_3)$ are related by $(K'_0,h',k')$, we have a homotopy commutative diagram
$$\xymatrix{K_1\ar[rr]^f&&K_2\ar[rr]^{f'}&&K_3\cr &K_0\ar[ul]^h\ar[ur]_k&&K'_0\ar[ul]^{h'}\ar[ur]_{k'}&\cr
&&K''_0\ar[ul]^{h''}\ar[ur]_{k''}&&}$$
where $(K''_0,h'',k'')$ is the pull-back of $\xymatrix{K_0 \ar[r]^{k}& K_2 &K'_0\ar[l]_{h'}}$. This is a cell complex (the common subdivision
of $K_0$ and $K'_0$), and $h''\colon K''_0\to K_0$ is a cellular homeomorphism, so $(K_1,c_1)$ and $(K_3,c_3)$
are related by $(K''_0, h\circ h'', k'\circ k'')$.
We choose any geometric basis for $(K''_0,p''_0)$, where $p''_0$ is the pull-back of the covering
$c_2^*(p)$ by $k\circ h'' \simeq k'\circ h''$. It follows that
$$\newt(\xymatrix{K_1 \ar@{~>}[r] & K_3})
=\newt(\xymatrix{K_1 \ar@{~>}[r] & K_2}) + \newt(\xymatrix{K_2 \ar@{~>}[r] & K_3})$$
by the composition
formula for the absolute torsion. Therefore
$(K_1,c_1)\sim (K_2,c_2)$ and $(K_2,c_2)\sim (K_3,c_3)$ implies
$(K_1,c_1)\sim (K_3,c_3)$.
\end{proof}
\begin{definition}
Let $\cS_\bullet (X,p)$ denote the set of equivalence classes of pointed torsion structures on $(X,p)$ for any polyhedron $X$.
We will use the notation  $[K,c]$ (or just $[K]$ when the reference map is understood) for the equivalence class of a torsion
structure $(K,c)$. If the covering space $p$ is understood, we will use the notation $\cS_\bullet (X)$ for short.
\end{definition}
 The basic idea of the equivalence relation is to identify two
pointed torsion structures whenever the homeomorphism  $c_2^{-1}\circ c_1\colon K_1 \to K_2$ is homotopic to a \emph{cellular} homeomorphism with zero absolute torsion. The more general formulation above will be useful in dealing with subdivisions or amalgamations of cell complexes. We remark that if $K \subset L$ is a subcomplex (containing the base-point), then a pointed torsion structure on $L$ induces a pointed torsion structure on $K$ by restriction of the data.

We now define a pairing
$$\Phi\colon \cS_\bullet (X,p)\times \cS_\bullet (X,p) \to K_1(\bZ[\pi])$$
by the formula
$$\Phi([K], [L]) = \newt(\xymatrix{K \ar@{~>}[r] & L})$$
where $(K,c_K)$ and $(L, c_L)$ are pointed torsion structures on
$(X,p)$ representing the classes $[K]$ and $[L]$ respectively.
This is well-defined, and may be computed by taking the absolute torsion of a cellular
approximation to $f= c_L^{-1}\circ c_K$ on a subdivision of $K$. Here are some basic properties of the pairing.
\begin{lemma} For any pointed torsion structures $[K]$, $[L]$ and
$[N]$ on $(X,p)$:
\begin{enumerate}
\item $\Phi([K], [N]) = \Phi([K], [L])+ \Phi([L], [N])$
\item $\Phi([K], [L]) = -\Phi([L], [K])$
\item $[K]$ is equivalent to $[L]$ if and only if $\Phi([K],[L])=0$.
\end{enumerate}
\end{lemma}
\begin{proof} The formulas follow directly from the definitions. The details are left to the reader.
\end{proof}
\section{ Absolute torsion structures on manifolds\label{four}}

In Theorem \ref{manifold} we give a formula for the absolute torsion of a closed, oriented manifold in terms of the $\Phi$-pairing. This formula will be used in the proof of Theorem B.
Our starting point is the fact that any closed $PL$ manifold $X$
of dimension $n$ has
a ``normal" triangulation (see \cite[\S 68]{seifert-threlfall}). More precisely,
 any simplicial $n$-complex $K$ homeomorphic to $X$ has the following properties:
\begin{enumerate}
\item Each $k$-simplex ($k <n$) of $K$  is incident with at least one
$n$-simplex (this is called a \emph{pure} $n$-complex).
\item Each $(n-1)$-simplex of $K$ is incident with exactly two $n$-simplices.
\item Any two $n$-simplices of $K$ can be connected by a sequence of
successively incident simplices of dimensions $n$ and $n-1$.
\end{enumerate}
By subdivision, we can find such a normal triangulation of arbitrarily
small diameter. If $X$ is oriented, we can choose a consistent
orientation for each $n$-simplex of $K$.  In addition, for any element
of $\cS_\bullet(X,p)$ we may find a representative $[K]$ such that the
underlying triangulation is normal. There is an analogous cell complex version of this notion, which we call a \emph{normal cell structure}
on $X$.

To compute the absolute torsion of $X$, we will use the  classical approach to Poincar\'e duality
(see \cite[Theorem 2.1]{wall-book}, \cite[\S~69]{seifert-threlfall}).
Let $K^\ast$ denote  the \emph{dual cell} complex
 associated to $K$. This is defined in terms of the triangulation, starting with the first barycentric subdivision $K'$ of
 $K$. Let $\Delta_r(K)$ denote the set of $r$-simplices of $K$.
 The vertices of $K'$ are the barycentres $\hat\sigma^r$
 of the simplices $\sigma^r \in \Delta_r(K)$, and its simplices
 have the form $[\hat\sigma^{r_0}, \hat\sigma^{r_1}, \dots,
 \hat\sigma^{r_s}]$, where  $\sigma^{r_j}$
 is a face of $\sigma^{r_{j+1}}$ for each $j$.
 Let $\Delta_r(K')$ denote the set of $r$-simplices of $K'$.
 The simplex $\sigma^r$
 of $K$ is the union of the simplices in $\Delta_r(K')$ which
 terminate with $\hat\sigma^{r}$.
 The cells $\cD \sigma^r$ of $K^\ast$ correspond bijectively to the
 simplices of $K$, and  $\cD \sigma^r$ is the union of the
 simplices of $K'$ which begin with $\hat\sigma^r$. Then
$ \dim \cD \sigma^r = n-r$ and $\sigma^r$ intersects $\cD\sigma^r$ transversely only at $\hat\sigma^r$. We denote
the set of $(n-r)$-cells of $K^\ast$ by $\Delta_{n-r}(K^\ast)$.

For a suitable chain approximation
to the diagonal in $K \times K^\ast$, cap product with the
fundamental cycle  takes the cochain on $K$ dual to
$\sigma^r$ to the chain $\cD \sigma^r$ on $K^\ast$.

\begin{definition}
\label{definition:dual_structure} Let $(K,c)$ be a torsion
structure on $X$. We define the \emph{dual torsion structure}
$(\cD K,\cD c)$ to be the torsion structure on $X$ with underlying
cell-complex $\cD K$ described above and with geometric basis as
follows:
\begin{enumerate}
\item The ordering of the cells on $\cD K$ is induced from the ordering
  of the cells on $K$:  for cells $\cD\sigma^r$ and $\cD\tau^r$ of the
  same dimension, $\cD\sigma^r$ comes before $\cD\tau^r$ if and only
  if $\sigma$ come before $\tau$.
\item The lift of each cell of $\cD K$ is induced by the lift of the
  corresponding dual cell.
  \item The cells of $\cD K$ are oriented such that cap product  with the
  fundamental cycle induces the based identity map
\[\cD\varphi_0\colon C(K)^{n-r} \to C(\cD K)_r\]
\end{enumerate}
With the above choice of geometric basis the chain complex of $\cD K$
is $C(\cD K) = C(K)^{n-*}$.  We choose the sign $\eta_{\cD K}$ to be
$\eta_{\cD K} = \eta_{C(K)^{n-*}}$.  Hence $C(\cD K) = C(K)^{n-*}$ as
  signed complexes (identified by $\cD\varphi_0$).
\end{definition}

\begin{theorem}\label{manifold} Let $X$ be a closed, oriented manifold of dimension $n$. Then
$$\newt(X) = \Phi([\cD K],[K]) \in \widehat H^n(\cy{2};K_1(\bZ[\pi_1(X,x_0)]))$$
\end{theorem}
\begin{proof}
We first choose a geometric basis for the first barycentric
subdivision $K'$ of $K$.  Since $X$ is an oriented manifold we may
apply the symmetric construction of Ranicki \cite{ra10}, \cite{ra11} to
find a symmetric Poincar\'e complex $(C(K),\varphi)$;  in particular, the
map $\varphi_0\colon  C(K)^{n-*} \to C(K)$ is a chain level representative of
the Poincar\'e duality map on $X$ given by cup-product with the
fundamental class.  We have a diagram of signed chain complexes which
commutes up to chain homotopy:
$${\vcenter
{\xymatrix
{C^{n-*}(K)\ar@{=}[r]^{\cD\varphi_0} \ar[d]_{\varphi_0}&
 C(\cD K)\ar[d] \\
C(K)\ar[r] & C(K')}
}}$$
Hence
\begin{eqnarray*}
\newt(X)&=&\newt(\varphi_0)\\
&=&\newt(C(\cD K)\to C(K'))-\newt(C(K) \to C(K'))\\
&=&\Phi([\cD K],[K']) - \Phi([K],[K'])\\
&=&\Phi([\cD K],[K])\ .
\end{eqnarray*}
\end{proof}

\section{ Fibre transport and transfer\label{transport}}

In this section we will describe a suitable algebraic setting for the fibre transport in a fibre bundle, and establish some useful properties for use in Section \ref{six}.
Let $\bA$ and $\bB$ be (small) additive categories with involution, and recall that
$\SPD_n(\bB)$ denotes the full sub-category of $\SD(\bB)$ consisting of signed $n$-dimensional chain complexes $C$  in $\bB$, such that
$C^{n-*}$ is chain equivalent to $C$, and $\chi(C) =0$ if $n$ is odd.
This sub-category inherits an involution via $C \mapsto C^{n-*}$ (see Definition \ref{fivefifteen}).

\begin{definition} A \emph{transfer functor} is a set of additive functors $J^D \colon \bA \to \SPD_n(\bB)$, indexed by a set ${_{\bA}Mod_{\bB}}$ of  objects $D$
in $\SPD_n(\bB)$ called $\bA$-$\bB$ \emph{bimodules}, such that
\begin{enumerate}
\item If $ D {\in}\, _{\bA}Mod_{\bB}$ and $D\cong D'$, then
$ D' {\in}\, _{\bA}Mod_{\bB}$.\\
\item For each isomorphism $g\colon D \to D'$ in $\SPD_n(\bB)$ there is a natural transformation $g_*\colon J^D \to J^{D'}$, with
$(id)_* = id$ and $(h\circ g)_* = h_*\circ g_*$ for any isomorphism $h\colon D' \to D''$. \\
\item For each object $M \in \bA$, $J^D(M)^* \cong J^{D^{n-*}}(M^*)$.
\hfill\qed
\end{enumerate}
\end{definition}
If $(C, \eta_C)$ is a signed complex in $\bA$, then $J^D(C)$ is a signed chain complex in $\SPD_n(\bB)$, with $\eta_{J^D(C)} = J^D(\eta_C)$. This formula gives an allowable sign because
$$J^D(\epsilon(M,N)) = \epsilon(J^D(M), J^D(N))$$ since $J^D$ is additive and
$\epsilon(M,N) = \tauiso(M\oplus N \to N \oplus M)$ for any objects $M$, $N$ of $\bA$. The composition
$$\trf_D=i_*\circ J^D_*\colon K_1^{iso}(\bA) \to K_1^{iso}(\bB)$$
is called the \emph{transfer} induced by $J^D$.

The main examples for our applications arise from the fibre transport.
\begin{example}[Fibre Transport] Let $R$ and $S$ be rings with involution, and let $D$ be an object in $\SPD_n(\bA(S))$. Suppose that
$$p\colon R \to [D,D]^{op}$$
is a homomorphism of rings with involution, where $[D,D]^{op}$
denotes the set of chain homotopy classes of self-chain maps
$D \to D$ over $S$. Then we let
$_{\bA(R)}Mod_{\bA(S)}$ be the collection of all objects $D' \cong D$ . Define
$$J^D \colon \bA(R) \to \SPD_n(\bA(S))$$
on objects by setting
$J^D(R) = (D, \eta_D)$, and extending additively. As a chain complex, $J^D(M) = M\otimes_R D$ for each object $M$ in $\bA(R)$, and the sign $\eta_{J^D(M)}$ is determined by the usual formula for direct sums.
If $f\colon M \to M'$ is a morphism in $\bA(R)$, then we get an induced morphism
$$f\otimes 1\colon M\otimes_R  D \to M'\otimes_R D,$$ where $p$ is used to define a left $R$-module structure on $D$.
If $g\colon D \to D'$ is an isomorphism in $\SPD_n(\bA(S))$, the natural transformation $g_*\colon J^D(C) \to J^{D'}(C)$ is defined by
$$1\otimes g\colon C \otimes_R D \to C \otimes_R D'\ .$$
The duality property (iii) follows from the identification $M = R^m$ as a based, free $R$-module. The involution $M \mapsto M^*$ maps the given base to the dual base. Under these identifications
$$M^*\otimes_R D^{n-*} = \oplus_{i=1}^m D^{n-*} = (\oplus_{i=1}^m D)^* = (M\otimes_R D)^*$$
In this example, the transfer induced by $J^D$ is usually denoted
$p^!\colon K_1(R) \to K_1(S)$ (see\cite{lueck1}).
\end{example}
As a special case (when $p$ is trivial) we get a transfer functor
$$J^D \colon \bA(R) \to \SPD_n(\bA(R\otimes S))$$ for
any two rings with involution (tensor product over $\bZ$).
\begin{example}[Tensor products] Let $C$ and $D$ be signed chain complexes over $\bA(R)$ and $\bA(S)$ respectively, and let
 $F_*(C\otimes D)$ denote the filtered complex structure defined in
Example \ref{tp_one}, Section \ref{three}. Then the \emph{signed tensor product}  is the signed filtered complex whose underlying chain complex is $C\otimes D$ and whose associated graded complex is $J^D(C, \eta_C) = (\bG_*(C\otimes D), J^D(\eta_C))$. The `internal' sign $\eta_{G_r(C\otimes D)}$ is the sign of the chain complex $J^D(C_r)$ over $S$, as defined above. In this example, the transfer $p^!= \trf_D$ is just multiplication by $\chi(D)$.
\hfill\qed
\end{example}

The functorial properties of $J^D$ give two useful formulas.
\begin{lemma}\label{trans_formulas}
Let $J^D\colon \bA \to \SD_n(\bB)$ be a transfer functor.
\begin{enumerate}
\item For each chain equivalence $f\colon C \to C'$ in $\bA$,
$$i_*\newt(J^D(f)) =i_*J^D_*(\newt(f)) =  \trf_D(\newt(f)) \in K_1^{iso}(\bB)\ .$$
\item If $J^D$ is induced by a fibre transport $R \xrightarrow{p} [D,D]^{op}$, then
for each chain equivalence $g\colon D \to D'$ in $\bA(S)$,
$$i_*\newt( J^{D}(C) \xrightarrow{g_*} J^{D'}(C)) = \chi(C)\cdot \newt(g) \in K_1^{iso}(\bA(S))\ .$$
\end{enumerate}
\end{lemma}
\begin{proof} Part (i)
 follows directly from the definitions. For Part (ii) observe that
 $g_*\colon J^D(C)\to J^{D'}(C)$ is an isomorphism of
  chain complexes in $\SD(\bA(S))$ and hence
\begin{eqnarray*}
&&i_*\newt(J^{D}(C) \xrightarrow{g_*} J^{D'}(C))\\[1 ex]
&&\hspace{1cm}\begin{array}{ll}
=&i_*\left(\sum(-)^r\tau^{iso}(1\otimes
g\colon C_r\otimes_R D
\to C_r\otimes_R D')\right) \\[1 ex]
=&\sum(-)^r\newt(\bigoplus_{\mathrm{rank}_R(C_r)}(g\colon D
\to D'))\\[1 ex]
=&\chi(C)\newt(g\colon D \to D') \in K_1^{iso}(\bA(S))
\end{array}
\end{eqnarray*}
\end{proof}

We now prove the product formula for the absolute torsion of symmetric
Poincar\'e structures announced in \cite{ak1}.
\begin{proposition}\label{product}
Let $(C,\varphi^C)$ and $(D,\varphi^D)$ be symmetric Poincar\'e complexes over
rings with involution $R$ and $S$ respectively.   Then
$$\newt(C\otimes D,\varphi^C\otimes
\varphi^D)=\chi(C)\newt(D,\varphi^D)+\chi(D)\newt(C,\varphi^C)$$
evaluated in $K_1(R\otimes S)$.
\end{proposition}
\begin{proof}
Let $\dim C =k$ and $\dim D = n$.
By definition
$$\newt(C\otimes D,\varphi^C\otimes \varphi^D) =
\newt(\varphi_0\otimes \phi_0\colon (C\otimes D)^{n+k-*} \to (C\otimes
D))\ .$$  The map $\varphi^C_0\otimes \varphi^D_0$ is given by the composition:
\[(C\otimes D)^{n+k-*}\xrightarrow{\theta_{C\otimes D}} C^{k-*}\otimes
D^{n-*} \xrightarrow{\varphi_0\otimes 1} C\otimes D^{n-*}
\xrightarrow{1\otimes\varphi_0}C \otimes D\]
The result now follows from Lemma \ref{trans_formulas}
 and
Lemma \ref{product_dual}.
\end{proof}

\section{ Absolute torsion structures on fibre bundles\label{five}}

The notion of pointed torsion structures will now be extended to
PL or smooth fibre bundles $F^n \xrightarrow{q} E^{n+k} \xrightarrow{p} B^k$ (see \cite[p.~181]{anderson1} for the definition of a PL fibre bundle). We always
assume that $F$, $E$ and $B$ are compatibly oriented PL or smooth, closed
manifolds (considered as polyhedra via the canonical $PL$ structure
compatible with their smooth structure). In both cases there exists a triangulation of the fibre bundle compatible with a given triangulation on the base (see  \cite{putz} for the smooth case). The orientation assumption means that $\pi_1(B,b_0)$ acts trivially on $H_n(F;\bZ)$ by fibre transport, and the isomorphism
$$H_{n+k}(E,\bZ) = H_k(B;\bZ) \otimes H_n(F;\bZ)$$
is compatible with the given orientations on total space, base and fibre.

  Our goal is to define a
notion of a \emph{pointed fibre bundle torsion structure} (PFBTS)
which will have as data a PTS on both the base and fibre, and which
will determine a PTS on the total space $E$.

Fix a base-point $b_0 \in B$ and a base-point
$e_0 \in E$ with $p(e_0)=b_0$. We will always assume that $B$ is connected.
In the definitions below, the geometric bases for torsion structures on $B$ or $E$ will use their \emph{universal} covering spaces (suppressed in the notation, together with the chosen lifts of the base-points).
In particular, if $c\colon K \to B$ is a pointed torsion structure (PTS), and
$\gamma \subset \bd \beta$ for some cell $\beta\in K$, then we have a component of the boundary chain map
$$\bd_{\beta, \gamma}\colon C(\widetilde \beta) \to SC(\widetilde \gamma)$$
which can be identified with the module homomorphism
$$d_{\beta, \gamma}\colon \bZ[\pi_1(B,b_0)] \to \bZ[\pi_1(B,b_0)]$$
given by multiplication with some element $d_{\beta, \gamma}
\in \{\pm\pi^{ab}\}$. This is a module isomorphism. We will use the images
$$p^!(\tau(d_{\beta, \gamma})) \in K_1(\bZ[\pi_1(E, e_0)])$$
under the fibre bundle transfer $p^!$ in the discussion below.

Let $K \to B$ and $L \to E$ be cell structures on the base and total
spaces such that the induced map $L \to K$ is cellular.  We call this a \emph{cell structure} on the fibre bundle. Then we may
regard $C(\widetilde{L})$ as a filtered complex by defining:
$$F_rC(\widetilde{L}) = C(p^{-1}(K^{\langle \leqslant r \rangle}))$$
where $K^{\langle \leqslant r \rangle}$ is the subcomplex of $K$
consisting of the cells of dimension $\leqslant r$. If $\dim F =
n$ and $\dim B = k$, then we obtain an admissible
$(n+k)$-dimensional $k$-filtered chain complex as defined in
Section \ref{three}.
 The associated complex of $F_*C(\widetilde{L})$
has the direct sum decomposition:
$$G_rC(\widetilde{L}) = S^{-r}\bigoplus_{\beta \in
  K^{(r)}}C(p^{-1}(\beta))$$
  indexed over the cells in $K$ of dimension $r$.
We will use the notation $C(\widetilde{F_\beta}) : = S^{-r}C(p^{-1}(\beta))$ for
  $\beta \in K^{(r)}$, so we have the identification
$$G_rC(\widetilde{L}) =\bigoplus_{\beta \in
  K^{(r)}}C(\widetilde{F_\beta})$$
Note that the cell structure on $E$ gives a cell structure on
  $F_0:=p^{-1}(b_0)$.
  Let $[F]$ denote a pointed torsion structure for  $F_0$, with respect to the \emph{pull-back} of the universal covering of $E$ by $q\colon F \to E$.
Fibre transport defines a ring morphism
\[\bZ[\pi_1B] \to [C(\widetilde{F_0}),C(\widetilde{F_0})]^{op}\]
so we have a transfer functor
\[J^{C(\widetilde{F_0})}\colon \bA(\bZ[\pi_1B]) \to \SPD_n(\bA(\bZ[\pi_1E]))\]
To shorten the notation, we will let $J^F:= J^{C(\widetilde {F_0})}$ when the PTS on $F$ is clear from the context.

  We can apply  this additive functor to the identification
  $$C_r(\widetilde K) = \bigoplus_{\beta \in
  K^{(r)}}C(\widetilde\beta)$$
  to obtain
  $$J^{C(\widetilde {F_0})}(C_r(\widetilde K)) = \bigoplus_{\beta \in
  K^{(r)}}J^{C(\widetilde {F_0})}(C(\widetilde\beta)) =
  \bigoplus_{\beta \in
  K^{(r)}} C(\widetilde {F_0})$$
  One of our main goals to construct a chain equivalence
\[\cE\colon \bG_*(C(\widetilde{L})) \to J^{C(\widetilde{F_0})}C(\widetilde{K})\]
by constructing chain equivalences
\[\cE^\beta\colon C(\widetilde{F_\beta}) \to C(\widetilde{F_0})\]
for each $\beta \in K$,  and defining $\cE$ as the direct sum of the $\cE^\beta$ over all such $\beta$.  This can
be achieved if we can construct such maps which are compatible with
boundary inclusions $\gamma < \beta$.  We will first construct the maps $\cE^\beta$
 for a special type of PTS on the base $B$ which comes from a
normal unfolding (defined below).  This construction will then be
extended to any PTS on the base $B$.

A pointed torsion structure $c_K\colon K \to B$ is \emph{small} with respect to the bundle $E\xrightarrow{p} B$, if $c(\beta)$ is contained in a locally trivial bundle neighbourhood for each $\beta \in K$.
\begin{definition}
A \emph{geometric basis} for a pointed fibre bundle
$(E,p, e_0)$ consists of the following data:
\begin{enumerate}
\item a pointed torsion structure $(K, c_K, k_0)$ on $B$, which is small with respect to the bundle $(E,p)$.
\item a cell structure $ L \to K$ on  $E \to B$, over  $c_K\colon K \to B$.
\item a geometric basis for  $F_\beta \subset L$, for each $\beta \in K$, with respect to the pull-back of the universal covering of $E$.
\hfill\qed
\end{enumerate}
\end{definition}
Notice that a geometric basis for $(E,p)$ and the fixed base-points give a PTS for the fibre $F$, represented by $(F_0, c_L|_{F_0})$. We denote this pointed torsion structure by $[F] \in \cS_{\bullet}(F)$.
\begin{definition}
A \emph{pointed fibre bundle torsion structure} (PFBTS)
on $(E,p)$
consists of a geometric basis $L \to K$ for  $(E,p, e_0))$, such that:
\begin{enumerate}
\item for each $\beta \in K$, a chain equivalence
$\cE^\beta\colon C(\widetilde{F_{\beta}}) \to C(\widetilde{F_0})$ where the map $\cE^{b_0}\colon C(\widetilde{F_0}) \to C(\widetilde{F_0})$ is the identity.\\
\item for each $\beta \in K$,
$$\newt(\cE^\beta\colon C(\widetilde{F_{\beta}}) \to C(\widetilde{F_0})) = 0 \in K_1(\bZ[\pi_1(E,e_0)])$$
\item the sum
\[\cE=\bigoplus\limits_{\beta \in K}\cE^{\beta}\colon \bG_*(C(\widetilde{L}))
\to J^{C(\widetilde{F_0})}C(\widetilde{K})\]
is a chain isomorphism in $\SD(\bA(\bZ[\pi_1(E,e_0)]))$.\hfill\qed
\end{enumerate}
\end{definition}
We will denote a PFBTS for $(E,p)$ by $(L\to K, \cE)$, and the induced PTS on the fibre by $[F]$. We will say that the PFBTS \emph{realizes} the given  structures  on base and fibre.

A PFBTS on $E\to B$ determines a PTS on $E$ as follows.  The choice of geometric basis on each fibre $F_\beta$ determines lifts and orientations of each of the
cell in $L$.  The associated complex $\bG_*(C(\widetilde{L}))$ may be
made into a signed complex by setting
\[G_r(C(\widetilde{L})) = \bigoplus_{\beta \in K^{\langle r
    \rangle}}C(\widetilde{F_\beta})\]
as signed complexes, and setting $\eta_{\bG_*C(\widetilde{L})} =
J^{F}_*(\eta_K)$.  We then give $C(\widetilde{L})$ the
filtered sign determined by the signs for $\bG_*(C(\widetilde{L}))$ and
define $\eta_L = \eta_{F_*C(\widetilde{L})}$.

\begin{lemma} $i_*\newt(\cE_*\colon \bG_*(C(\widetilde{L})) \to J^{C(\widetilde{F_0})}C(\widetilde{K}))= 0 \in K_1(\bZ[\pi_1(E,e_0)])$.
\end{lemma}
\begin{proof} This follows directly from the sign conventions chosen, and the formula
$$\newt(\cE_*) = \sum_r (-1)^r\tauiso(\cE_r\colon G_r(C(\widetilde L)) \to J^F C_r(\widetilde{K})) + \eta_{J^F C(\widetilde{K})} -
\eta_{\bG_*(C(\widetilde{L}))}$$
for the absolute torsion of a chain isomorphism.
\end{proof}

\section{ Absolute torsion of fibre bundles\label{six}}

We will show in Theorem \ref{realization} that every  smooth (or PL) fibre bundle $F \to E \to B$ as above admits a pointed fibre bundle pointed torsion structure.
First we define the \emph{normal unfolding} of a normal triangulation
of a closed  manifold. The idea is to cut $X$ open along
certain $(n-1)$-simplices of the triangulation $K$ to obtain a contractible
complex projecting onto $X$.
\begin{definition}
 Let $K=\bigcup \{ \Delta^n_\alpha \vv \alpha \in J_K\}$, where the $\Delta^n_\alpha$
are the distinct $n$-simplices in the triangulation of $X$.
Let $j_\alpha\colon \Delta^n_\alpha \to K$ denote the inclusion
maps of the $n$-simplices.
 A \emph{normal unfolding} of $K$ is a pure simplicial $n$-complex
$\wK$ equipped with a simplicial map $f\colon \wK \to K$, such that
\begin{enumerate}
\item $\wK$ is a union of $n$-simplices, with index set $J_{\wK}$, and each $(n-1)$-simplex of $\wK$
is incident with one or two $n$-simplices.
\item There is a bijection $\theta\colon J_{\wK} \to J_K$ such that
$f|\Delta^n_\alpha = j_{\theta(\alpha)}$ for all $\alpha \in J_{\wK}$.
\item $\wK$ is contractible.
\end{enumerate}
\end{definition}
We again have an analogous notion for normal cell structures.
\begin{lemma} Let $X$ be a closed smooth manifold of dimension
$n$. Each normal cell structure $K$ on $X$ admits a
normal unfolding $(\wK, f)$.
\end{lemma}
\begin{proof} We give the proof for the case of triangulations. Consider the set of barycentres
$V=\{v_\alpha\vv \alpha \in J_K\}$ of the $n$-simplices
of $K$. We can join the barycentres of each pair $\{\Delta^n_\alpha, \Delta^n_\beta\}$
 of incident $n$-simplices with an edge $e_{\alpha\beta}$
 through their common $(n-1)$ face to obtain a connected graph (the dual 1-skeleton of $K$). Let $T = (V, E)$ denote a maximal tree in this graph, and note that $T$ has vertex set consisting of all the barycentres.
 Let $\wK$ be the quotient space of the disjoint union of the $n$-simplices in $K$, where
 we identify
 two $n$-simplices $\Delta^n_\alpha$,  $\Delta^n_\beta$ in $\wK$ along an
 $(n-1)$ face if and only if $e_{\alpha\beta} \in E$. By
 construction, $\wK$ is a thickening of the tree $T$,  so $\wK$ is
 contractible. The other properties of a normal unfolding
 are clear.
\end{proof}
A normal unfolding $\wK \to K$ can be used to make consistent choices of lifts to the universal covering for simplices in $K$. Any cell $\beta \in K$ is contained in a unique $n$-cell of minimal ordering. We can specify a unique lift $\hat\beta \in \wK$ by requiring $\hat\beta$ to be contained in the lift $\hat\Delta$ of $\Delta$.
\begin{definition}
A pointed  torsion structure $K \to B$ is \emph{normal} if for
some normal unfolding $\wK \to K$, the chosen lifts of adjacent
cells $\beta, \gamma \in  K$ are adjacent in the universal
covering $\widetilde K$ whenever their lifts $\hat\beta$ and
$\hat\gamma$ are adjacent in $\wK$.
\end{definition}

\begin{remark}
If $f\colon \wK \to K$ is a normal unfolding, the map $f$ induces
an equivalence relation on the $r$-simplices of $\wK$. Two $r$-simplices are equivalent
if they are identified by $f$ to the same simplex in $K$.
Let $\bd \wK$ denote the subset of $K$ consisting of
all the simplices whose equivalence class contains more than one element. This is a sub-complex of $\wK$ called its \emph{boundary}.

Let $v_0$ denote the barycentre of one of the $n$-simplices (chosen as a base point, so $f(v_0) = k_0$), and let $v_i \in \tau_i$ denote the barycentre of each of the $n$-simplices $\Delta^n_i$ with a face $\tau_i \subset \bd\wK$. Let $x_i\in \tau_i$ denote the barycentre of each such face, and $e_i = [v_i,x_i]\subset\Delta^n_i$ the linear path joining $v_i$ to $x_i$. Let $\sigma_i$ denote the unique path in the maximal tree joining $v_0$ with the barycentre  $v_i \in \Delta^n_i$, followed by the path $e_{i}$ to $x_i$. If $\tau_i$ and $\tau_j$, $i<j$, are  a pair of $(n-1)$-faces identified under $f\colon \wK \to K$, then $f(e_i) \cup f(e_j) = e_{ij}$ is an edge in the dual 1-skeleton of $K$.
The fundamental group $\pi_1(K, k_0)$ is generated by  loops of the form $u_{ij}: =f(\sigma_i)\cup e_{ij} \cup f(\sigma_j)$.
\hfill\qed
 \end{remark}
Here is our main result about the existence of pointed fibre bundle torsion structures.
\begin{theorem}[\textbf{Realization}]
\label{realization}
Let $E \to B$ be a fibre bundle of compatibly oriented closed PL manifolds.
Given a PTS on $K\to B$, a cell structure $L \to K$ on the fibre bundle, and a geometric basis
on $F$, there exists
 a PFBTS $(L \to K,\cE)$ realizing this data.
\end{theorem}
The main steps in the proof  are contained in  the following special case.
\begin{proposition}
\label{prop:special_pfbts}  There
 exists a PFBTS $(L \to K,\cE)$ on $E\to B$ realizing a given
normal pointed torsion structure $K \to B$  on $B$,
a given cell structure $L\to K$  on the fibre bundle,
 and  a given geometric basis for $F$.
\end{proposition}
\begin{proof}
We are required to construct chain maps $\cE^\beta\colon C(\widetilde{F_\beta})
 \to C(\widetilde{F_0})$ such that the direct sum $\cE\colon\bG_*(C(\widetilde{L}))
 \to J^FC(\widetilde{K})$ is a chain equivalence. Once these maps are constructed,  we can choose a geometric basis on each $F_\beta$
such that the maps $\cE^\beta\colon C(\widetilde{F_\beta}) \to
C(\widetilde{F_0})$ have
trivial absolute torsion (this is condition (ii) for a PFBTS).
 The condition that $\cE$ is a chain map
 is equivalent to the statement that for a pair of simplices $\gamma <
 \beta$ where $\gamma$ is an $(r-1)$-dimensional cell in the
 boundary of an $r$-dimensional cell $\beta$, the diagram
 \eqncount
\begin{equation}
\label{eqn:boundary}
\vcenter{
\xymatrix{
C(\widetilde{F_\beta}) \ar[d] \ar[r]
& J^{F}C(\beta) \ar[d] \\
C(\widetilde{F_\gamma}) \ar[r]
& J^{F}C(\gamma)
}}
\end{equation}
commutes.  We will first construct the maps $\cE$ and then show that
this diagram commutes for all $\gamma < \beta$.

Let $f\colon \wK \to K$ be the normal unfolding of $K$, and $g\colon
\wL \to L$ the pull-back of $L \to K $ over $\wK$.
The chain complex $\bG_*(C(\widehat{L}))$ splits over the cells
$\hat{\beta} \in \widehat{K}$ as:
\[\bG_*(C(\widehat{L})) = \bigoplus_{\hat{\beta}\in \widehat{K}}
C(\widetilde{F_{\hat{\beta}}})\]
We assume that
the base-point $k_0$ is uniquely covered ($f^{-1}(k_0) = v_0$).
The chain complex $C(\widehat{K})$ is a based $\bZ$-module chain
complex since we have an orientation for each cell.  Since
$\widehat{K}$ is contractible the pull-back of the universal cover of
$B$ over the map $\widehat{K}\to K\to B$ is trivial and hence we have a
chain map
\[f_*\colon C(\widehat{K})\otimes_\bZ \bZ[\pi_1B] \to C(\widetilde{K})\]
of based signed chain complexes over $\bZ[\pi_1 B]$.
On the total space we have a filtered chain map
$g_*\colon C(\widehat{L}) \to C(\widetilde{L})$
whose associated map
\[\bG_*(g_*)\colon \bG_*(C(\widehat{L})) \to \bG_*(C(\widetilde{L}))\]
splits over each cell $\beta$ of $K$,  and lifts $\hat{\beta}$ of
$\beta$ giving
an isomorphism of chain complexes:
\[g_{\hat{\beta}}\colon C(\widetilde{F_{\hat{\beta}}}) \to C(\widetilde{F_\beta})\]
Note that
\[J^F(C(\widehat{K})\otimes_\bZ \bZ[\pi_1B]) =
C(\widehat{K})\otimes_\bZ C(\widetilde{F_0})\ .\]
 Since
$\wK$ is contractible and the triangulation $K \to B$ is small with
respect to the bundle $E \to B$, there is a trivializing homeomorphism
$$\Psi_K\colon |\wL| \to |\wK| \times F$$
of the pull-back bundle $\hat p\colon \widehat L \to \widehat K$
which gives rise to a
chain homotopy commutative diagram of
chain maps:
\[\xymatrix{C(\hat{p}^{-1}(\hat{\beta})) \ar[d] \ar[r] &
  C(\widetilde{\beta})\otimes_\bZ
  C(\widetilde{F_0})=C(\widetilde{F_0})
  \ar[d] \\
C(\widehat{L})\  \ar[r] & \ \ \ C(\widehat{K})\otimes_\bZ C(\widetilde{F_0}) =
  J^F(C(\widehat{K})\otimes_\bZ \bZ[\pi_1B])
}\]
 for each $\hat{\beta}\in \widehat{K}$.
Let $r$ be the dimension of $\beta$.  By taking the associated
complexes and maps in this diagram, we
construct maps $\cE^{\hat{\beta}}$:
\[\xymatrix{C(F_\beta) \ar[d] \ar[r]^{\cE^{\hat{\beta}}} &
  C(\widetilde{F_0})   \ar[d] \\
G_r(C(\widehat{L})) \ar[r] & \ \ J^F(C_r(\widehat{K})\otimes_\bZ
\bZ[\pi_1B])
}\]
The sum $\widehat{\cE}\colon \bG_*(C(\widehat{L})) \to
J^{C(\widetilde{F_0})}(C(\widehat{K})\otimes_\bZ \bZ[\pi_1B])$ is a
chain map and each $\cE^{\hat{\beta}}$ is a chain equivalence.

Each cell $\beta$ in $K$ lies in an $n$-simplex of minimal ordering so
we have a particular lift $\hat{\beta}$ of $\beta$.  We define the
chain equivalence
$\cE_\beta$ to be
\[\cE^\beta = \cE^{\hat{\beta}}\circ (g_{\hat{\beta}})^{-1} \colon
C(\widetilde{F_\beta}) \to C(\widetilde{F_{\hat{\beta}}}) \to
C(\widetilde{F_0})\]
where $g_{\hat\beta}$ is the restriction of $g\colon \wL \to L$ to $\hat{\beta}$.
It remains to show that $\cE_*$ is a chain map, in other words that
the diagram  (\ref{eqn:boundary}) commutes for all pairs $\beta,\gamma$.
This is true if the lifts of the cells $\beta$ and
$\gamma$ are adjacent in $\widehat{K}$, since we have a commutative diagram:
\[\xymatrix{
C(\widetilde{F_\beta}) \ar[r]^{(g_{\hat{\beta}})^{-1}} \ar[d]&
C(\widetilde{F_{\hat{\beta}}})
\ar[r]^(0.35){\cE^{\hat{\beta}}}  \ar[d]
& \  C(\widetilde{F}) =J^FC(\beta) \ar[d]\\
C(\widetilde{F_\gamma}) \ar[r]^{(g_{\hat{\gamma}})^{-1}} &
C(\widetilde{F_{\hat{\gamma}}})
\ar[r]^(0.35){\cE^{\hat{\gamma}}} & \ C(\widetilde{F})=J^FC(\gamma)
}\]

 Suppose now that the lift of
$\gamma$ is not adjacent to that of $\beta$.  Then the component of
 the differential in $C(\widetilde{K})$ given by $C(\beta)
\to C(\gamma)$ is $d_{\beta,\gamma} =\pm t$ for some non-trivial group element $t  \in
\pi_1B$, when considered as a map of based complexes.  We will denote
by $\pm$ the sign in front of $t$.  Let
$\hat{\gamma}'$ be the lift of $\gamma$ which is adjacent
to the chosen lift of $\hat{\beta}$ of $\beta$.
%%%%%%%%%%%%%%%%%%
The two faces $\hat{\gamma}$, $\hat{\gamma}'$ are identified
under $\widehat{K} \to K$.  Let $f_{\hat{\gamma}}$,
$f_{\hat{\gamma}'}$  be paths in $\widehat{K}$ joining the
barycentres $x_{\hat{\gamma}}$, $x_{\hat{\gamma}'}$ with the
base point $v_0$.  Then the image of $f_{\hat{\gamma}}\cup
f_{\hat{\gamma}'}$ in $K$ is a loop representing $t$.  Hence fibre
transport around $t$ is described by the following composition:
\[F \to F\times{x_{\hat{\gamma}}}
\xrightarrow{\Psi|_{x_{\hat{\gamma}}}}
\hat{p}^{-1}(x_{\hat{\gamma}}) \to
p^{-1}(x_{\gamma}) \to
\hat{p}^{-1}(x_{\hat{\gamma}'})
\xrightarrow{(\Psi|_{x_{\hat{\gamma}'}})^{-1}}
F\times{x_{\hat{\gamma}'}} \to F\]
where $x_\gamma$ is the barycentre of $\gamma$ in $K$.
%%%%%%%%%%%%%%%%%%%%
Therefore the composition of homotopy equivalences:
\[F \to F\times |\hat{\gamma}|
\xrightarrow{\Psi_K|_{\hat{\gamma}}} F_{\hat{\gamma}}
\xrightarrow{g|_\gamma}
F_{\gamma} \xrightarrow{(g|_{\gamma'})^{-1}} F_{\hat{\gamma}'}
  \xrightarrow{(\Psi_K|_{\hat{\gamma}' })^{-1}} F\times
  |\hat{\gamma}'| \to F
\]
is described up to homotopy by fibre transport around $t$.  Hence the
composition
\[C(\widetilde{F}) \xrightarrow{(\cE^{\hat{\gamma}})^{-1}}
C(\widetilde{F_{\hat{\gamma}}})
\xrightarrow{g_{\hat{\gamma}}}  C(\widetilde{F_{\gamma}})
\xrightarrow{(g_{\hat{\gamma}'})^{-1}}
C(\widetilde{F_{\hat{\gamma}'}})
\xrightarrow{\cE^{\hat{\gamma}'}} C(\widetilde{F})\]
is given (up to chain homotopy) by $J^F(t)\colon C(\widetilde{F}) \to
C(\widetilde{F})$.  From this we deduce that
\[\cE^\gamma = J^F(t)\circ (
\cE^{\hat{\gamma}'} \circ (g_{\hat{\gamma}'})^{-1})\]
Since $\hat{\gamma}'$ is adjacent to $\hat{\beta}$ we know that
\[\xymatrix{
C(\widetilde{F_\beta}) \ar[d] \ar[r]^{\cE^\beta}
& C(\widetilde{F}) \ar[d]^{J^F(\pm)} \\
C(\widetilde{F_\gamma})\ \ \  \ar[r]^{\cE^{\hat{\gamma}'} \circ
  (g_{\hat{\gamma}'})^{-1}}
& \ \ C(\widetilde{F})
}\]
commutes.  It now follows that
 \[
 \xymatrix{
C(\widetilde{F_\beta}) \ar[d] \ar[r]^{\cE^\beta}
& C(\widetilde{F}) \ar[d]^{J^F(\pm t)} \\
C(\widetilde{F_\gamma}) \ar[r]^{\cE^\gamma}
& C(\widetilde{F})
}
\]
commutes as required.
\end{proof}

\begin{proof}[The proof of Theorem \ref{realization}]
Let $K' \to B$ be a normal PTS on $B$ with the same underlying
cell-structure
as $K$.  Then we can find a PFBTS $(L' \to
K',\cE')$ with $L'$ having the same underlying cell structure as $L$
and fibre structure $F$ by the above
proposition.  For each cell $\beta\in K$ we have a corresponding cell
$\beta'\in K'$ and an identification of $C(F_\beta)$ with
$C(F_{\beta'})$.  We use this identification to choose a geometric
basis for $F_0$ such that $C(F_{k_0})\to C(F_{k'_0})$ is the based
identity map.
We now define the maps $\cE^\beta$ by
\[\cE^\beta =
J^{F}(C(\beta')\to C(\beta)) \circ (\cE')^{\beta'}\colon C(\widetilde{F_\beta})
\to C(\widetilde{F_0})\]
using the identification of $C(F_\beta)$ with
$C(F_{\beta'})$.  We again choose a geometric basis on each $F_\beta$
such that the maps $\cE^\beta\colon C(\widetilde{F_\beta}) \to
C(\widetilde{F_0})$ have
trivial absolute torsion.  We claim that $(L \to K,
 \cE)$ is a PFBTS.  The only condition remaining to check is that
 $\cE$ is a chain map.  However, by construction, $\cE$ is the
 composition
\[\bG_*(C(\widetilde{L})) \to \bG_*(C(\widetilde{L'})) \xrightarrow{\cE'}
J^FC(\widetilde{K'}) \xrightarrow{J^F(C(\widetilde{K'})\to
  C(\widetilde{K}))} J^FC(\widetilde{K})\]
of chain maps and hence is itself a chain map.
\end{proof}
\begin{remark}
\label{remark:transport}
The maps $\cE^\beta$ in a PFBTS have a geometric interpretation.  For
an oriented cell $\beta \in K$, the intersection of cells of
$p^{-1}(\beta)$ with $p^{-1}(\hat{\beta})$, the barycentre of $\beta$,
defines a cell-structure on $p^{-1}(\hat{\beta})$.  A lift of a cell
in $p^{-1}(\beta)$ determines a lift of the corresponding cell in
$p^{-1}(\hat{\beta})$, and also the orientation of the cell in
$p^{-1}(\beta)$ along with the orientation of $\beta$ determines an
orientation of the cell in $p^{-1}(\hat{\beta})$.  This determines a
PTS on $p^{-1}(\hat{\beta})$; moreover we have an identification of
based chain complexes:
\[C(\widetilde{F_\beta}) \to C(p^{-1}(\hat{\beta}))\ .\]
We can regard the choice of lift of $\beta$ as a path $\sigma$ from the base
point $b_0$ to $\beta$.  Then the map $\cE^\beta$ composed with the
above identification 
\[\cE^\beta\colon C(p^{-1}(\hat{\beta})) \to C(\widetilde{F})\]
is a chain level representative of the map from $p^{-1}(\hat{\beta})$
to $F$ given by fibre transport along $\sigma$.
\end{remark}

\section{ The Proof of Theorem B\label{thmB}}

In order to prove Theorem B we will need to show how the absolute
torsions computed from a pointed torsion structure on a fibre bundle
vary when we change the PTS  on the base or on the fibre. Let
$\cS_\bullet(F \xrightarrow{q} E)$ denote the set of pointed torsion structures on $F_0 = p^{-1}(b_0)$ with respect to the pull-back of the universal covering of $E$.  We summarize the results as follows:

\begin{proposition}\label{propone}
Given  $[B] \in \cS_\bullet(B)$ and $[F] \in \cS_\bullet(F \xrightarrow{q} E)$, there exists a fibre bundle torsion structure $(L \to K, \cE)$ on $F\to E \to B$, such that $[K]
= [B]$ and the class of the fibre structure is $[F]$.  Moreover any
two such fibre bundle torsion structures  represent the same element
in
$\cS_\bullet(E)$.
\end{proposition}

Let $\fS_{[B][F]} \in \cS_\bullet(E)$  denote the induced torsion structure on the total space.

\begin{proposition}\label{proptwo}
Let $[B_1], [B_2] \in S_\bullet(B)$ be any two torsion structures for the base, and let
$[F_1], [F_2]\in \cS_\bullet(F \xrightarrow{q} E))$ be any two fibre structures.  Then
\[\Phi(\fS_{[B_1][F_1]},\fS_{[B_2][F_2]}) = p^!\Phi([B_1],[B_2]) +
  \chi(B)\Phi([F_1],[F_2])\]
\end{proposition}
The \emph{dual torsion structure} $\cD\fS_{[B][F]}$ of $\fS_{[B][F]}$ was defined in
Definition \ref{definition:dual_structure}.

\begin{proposition}\label{propthree}
Let $\cD\fS_{[B][F]}$ denote the dual torsion structure of $\fS_{[B][F]}$ in
$\cS_\bullet(E)$.  Then
\[\Phi(\cD \fS_{[B][F]},\fS_{[\cD B][\cD F]}) = 0\]
\end{proposition}
These properties will be established in a sequence of lemmas.

\begin{lemma}
\label{remark:unique_e}
Given PTSs on $K$ and $F$ the maps
$\cE^{\beta}\colon C(\widetilde{F_\beta})
\to C(\widetilde{F_0})$ are uniquely determined by the properties:
\begin{enumerate}
\item The map $\cE^{b_0}\colon C(\widetilde{F_0}) \to C(\widetilde{F_0})$ is the
  identity.
\item The direct sum $\cE = \bigoplus_\beta \cE^\beta$ is a chain equivalence.
\end{enumerate}
\end{lemma}
\begin{proof}
Suppose that $\{\cE_1^\beta\}$ and $\{\cE_2^\beta\}$ are two
sets of chain maps satisfying these conditions. Since $B$ is connected, any cell $\beta \in K$ can be connected to the base-point by a sequence of adjacent pairs of cells. If $\gamma \subset \bd \beta$, we have a commutative diagram
$$\xymatrix{C(\widetilde{F_\beta})\ar[r]^{\cE_1^\beta}\ar[d]& C(\widetilde{F_0})\ar@{=}[d]\\
C(\widetilde{F_\gamma})\ar[r]^{\cE_1^\gamma} &
C(\widetilde{F_0})
}
$$
and similarly for $\cE_2$. We are using the identification
$J^FC(\beta)=C(\widetilde{F_0})$, valid for every cell. But the horizontal maps are isomorphisms in $\SPD_n(\bZ[\pi_1(E,e_0)])$, so $\cE_1^\beta = \cE_2^\beta$ if and only if
 $\cE_1^\gamma = \cE_2^\gamma$. However, these two maps agree at the base-point $b_0$, so they agree over every cell and
 $\cE_1 = \cE_2$.
 \end{proof}
\begin{lemma} Suppose that  $(L_1 \to K, \cE_1)$ and $(L_2 \to K,
  \cE_2)$ are  pointed fibre bundle torsion structures on $
  E \to B$ with fibre structures $F_1$ and $F_2$ respectively. If
  $h\colon L_1\to L_2$ is a subdivision over $K$, then
  there is a  commutative diagram
$$\xymatrix{\bG_*(C(\widetilde L_1))
    \ar[d]_{\bG_*(h)}\ar[r]^{\cE_1}&J^{F_1}(C(\widetilde K))
\ar[d]^{h_*}\\
\bG_*(C(\widetilde L_2))\ar[r]_{\cE_2}&J^{F_2}(C(\widetilde K))}$$
in $\SPD_n(\bZ[\pi_1(E,e_0)])$.
\end{lemma}
\begin{proof} Since any subdivision over $K$ is a filtered map, we can pass to the associated graded complexes.
The composition
\[\cE'_1\colon \bG_*(C(\widetilde{L_1})) \xrightarrow{\bG_*(h)} \bG_*(C(\widetilde{L_2}))
\xrightarrow{\cE_2}
J^{F_2}C(\widetilde{K}) \xrightarrow{(h_*)^{-1}}
J^{F_1}C(\widetilde{K})\]
splits over each simplex $\beta\in K$ to the composition:
\[{\cE'_1}^\beta\colon C(\widetilde{F_{1\beta}}) \xrightarrow{\bG_*(h_*|_\beta)}
C(\widetilde{F_{2\beta}}) \xrightarrow{\cE_2^\beta} J^{F_2}C(\beta)
\xrightarrow{(h_*)^{-1}} J^{F_1}C(\beta)\]
In particular the map ${\cE'_1}^{b_0}$ coincides with
${\cE_1}^{b_0}$, so the above remark implies that $\cE'_1 = \cE_1$.
Hence the diagram commutes.
\end{proof}

\begin{corollary}\label{lemma: fibrechange}
\label{cor:fibre}
Suppose that $(L_1 \to K,\cE_1)$ and $(L_2 \to K,\cE_2)$ are
pointed fibre bundle torsion structures on a fibre bundle
$F\to E \to B$.
The induced pointed torsion
structures  $[L_1]$ and $[L_2]$ on $E$ satisfy
\[\Phi([L_1],[L_2]) = \chi(B)\Phi([F_1],[F_2])\]
\end{corollary}
\begin{proof}
If $L_2$ is a subdivision of $L_1$, then the commutativity of the
diagram implies that
\[\newt(\bG_*(C(\widetilde{L_1})) \to \bG_*(C(\widetilde{L_2}))) = \newt(
J^{F_1}(C(\widetilde{K})) \to J^{F_2}(C(\widetilde{K})))\]
since the maps $\cE_1$ and $\cE_2$ have trivial absolute torsion.
By Lemma \ref{trans_formulas}
\begin{eqnarray*}i_*\newt(\bG_*(C(\widetilde{L_1})) \to \bG_*(C(\widetilde{L_2})))
&=& \chi(B)\newt(C(\widetilde{F_1})\to C(\widetilde{F_2})) \\
&=& \chi(B)\Phi([F_1],[F_2])\ .
\end{eqnarray*} From Theorem \ref{invariance2} on the torsion of filtered maps we
see that:
\begin{eqnarray*}
\Phi([L_1],[L_2]) & = & \newt(C(\widetilde{L_1}) \to C(\widetilde{L_2}))\\
& = & i_*\newt(\bG_*(C(\widetilde{L_1})) \to \bG_*(C(\widetilde{L_2})))\\
& = & \chi(B)\Phi([F_1],[F_2])
\end{eqnarray*}
as required.  We now consider the general case where $L_2$ is not
necessarily a subdivision of $L_1$; in this case we may find a
filtered cellular map $L_3 \to L_1$ which is a common subdivision (as
cell-complexes) of $L_1$ and
$L_2$.  By Theorem \ref{realization} we may find a PFBTS whose
underlying triangulation is $L_3$, we will also denote this by $L_3$.
It now follows that
\begin{eqnarray*}
\Phi([L_1],[L_2]) &=& \Phi([L_2],[L_3]) - \Phi([L_1],[L_3]) \\
&=&\chi(B)\Phi([F_2],[F_3]) - \chi(B)\Phi([F_1],[F_3]) \\
&=&\chi(B)\Phi([F_1],[F_2])
\end{eqnarray*}
as required.
\end{proof}

\begin{lemma} Suppose that  $(L_1 \to K_1, \cE_1)$ and
$(L_2 \to K_2, \cE_2)$ are  pointed fibre bundle torsion
structures on $E
  \to B$. Let $h\colon L_1 \to L_2$, $f\colon K_1\to K_2$  be cellular
  homeomorphisms, with $p_{L_2}\circ h = f\circ p_{L_1}$ and $h$ inducing the identity on
   $F = p^{-1}(b_0)$. If
  both PFBTS have the same fibre structure under this
  identification, then there is a commutative
  diagram
$$\xymatrix{\bG_*(C(\widetilde L_1)) \ar[d]_{\bG_*(h)}\ar[r]^{\cE_1}&J^{F}(C(\widetilde K_1))
\ar[d]^{J^F(f)}\\
\bG_*(C(\widetilde L_2))\ar[r]_{\cE_2}&J^{F}(C(\widetilde K_2))}$$
in $\SD_n(\bZ[\pi_1(E,e_0)])$.
\end{lemma}
\begin{proof}
We first consider the case where $K_1$ and $K_2$ differ only by the
choice of PTS.  Then the composition
\[\bG_*(C(\widetilde{L_1})) \xrightarrow{\bG_*(h)} \bG_*(C(\widetilde{L_2}))
\xrightarrow{\cE_2} J^FC(\widetilde{K_2})
\xrightarrow{J^F(f_*^{-1})} J^FC(\widetilde{K_1})\] splits over
the cells of $K_1=K_2$ to give maps ${\cE'}_1^\beta\colon
\bG_*(C(\widetilde{F_\beta})) \to J^FC(\widetilde{K_1})$.  However
Remark \ref{remark:unique_e} shows that these maps must coincide
with the $\cE_1^\beta$, so the diagram commutes in this case.

It is now sufficient to consider the case where the map $K_1 \to K_2$
is an elementary subdivision of one cell $\beta$ and where $L_1 \to L_2$ is a
cellular isomorphism away from $p^{-1}(\beta)$.  We may also assume
 the geometric basis for $K_2$ is the same as that for $K_1$ away
from $\beta$ and that the geometric basis for each $C(\widetilde{F_{1\gamma}})$
coincides with that for $C(\widetilde{F_{2\gamma}})$ for $\gamma \neq \beta$.
We may also assume that the
maps $\cE_1^\gamma$ and $\cE_2^\gamma$ are equal for all $\gamma\neq\beta$
so it is sufficient to
show that the diagram commutes on the $C(\widetilde{F_\beta})$ factor.
Let $\{\beta_i\}_i$ be the components of the subdivision of $\beta$ in
$K_2$ which have the same dimension as $\beta$.  Then we must show
that the diagram
\[\xymatrix{C(\widetilde{F_\beta}) \ar[d]_{\bG_*(h|_\beta)} \ar[r]^{\cE^\beta} &
  J^FC(\beta) \ar[d]^{J^F(f|_\beta)}
  \\
  \ \bigoplus_iC(\widetilde{F_{\beta_i}}) \ar[r]^{\oplus\, \cE^{\beta_i}}
  &\  \bigoplus_iJ^FC(\beta_i)}
\]
commutes.  The RHS of the diagram splits over the components
$C(\widetilde{F_{\sigma_i}})$ so it is sufficient to show that
\[\xymatrix{C(\widetilde{F_\beta}) \ar[d]_{\bG_*(h|_{\beta_i})} \ar[r]^{\cE^\beta} & J^FC(\beta) \ar[d]^{J^F(f|_{\beta_i})}
  \\ C(\widetilde{F_{\beta_i}}) \ar[r]^{\cE^{\beta_i}} & J^FC(\beta_i)}
\]
commutes for each $i$.  Let $\tau$ be a boundary component of both
$\beta_i$ and $\beta$ (such a component can always be found for an
elementary subdivision).  Then we have a diagram of chain equivalences:
\[\xymatrix{C(\widetilde{F_{\tau}})
 \ar[r]\ar[d]^{\cE^{\tau}} &
\ C(\widetilde{F_\beta})\ \ar[r]^{\bG_*(h|_{\beta_i})}\ar[d]^{\cE^\beta}&
\ C(\widetilde{F_{\beta_i}})\ \ar[r]\ar[d]^{\cE^{\beta_i}}&
C(\widetilde{F_{\tau}})\ar[d]^{\cE^{\tau}}\\
  J^FC(\tau) \ar[r]  &
\ J^FC(\beta)\ \ar[r]^{J^F(f|_{\beta_i})} &
\ J^FC(\beta_i)\ \ar[r]&
J^FC(\tau)
}\]
The outer square commutes since the compositions along the top and bottom are
the identity maps.  The left and right squares commute since both
$\cE_1$ and $\cE_2$ are chain maps.  Hence the middle square commutes
as required.
\end{proof}

\begin{corollary}
\label{cor:base}
\label{bundle_structure_base_change}
Let $K_1$ and $K_2$ be cell-decompositions of $B$; let $L_1 \to K_1$
and $L_2 \to K_2$ be PFBTSs on $E \to B$ with fibre
class $[F]$.  Then the induced PTSs on $E$, $[L_1]$
and $[L_2]$ satisfy
\[\Phi([L_1],[L_2]) = p^!\Phi([K_1],[K_2])\]
\end{corollary}
\begin{proof}
If $L_2$ and $K_2$ are subdivisions of $L_1$ and $K_1$ respectively,
then commutativity of the diagram of the above lemma and Lemma
\ref{trans_formulas} give:
\begin{eqnarray*}
i_*\newt(\bG_*(C(\widetilde{L_1}))\to \bG_*(C(\widetilde{L_2}))) &=&
\newt(J^FC(\widetilde{K_1})\to J^FC(\widetilde{K_2})) \\
&=& p^!\newt(C(\widetilde{K_1}) \to C(\widetilde{K_2}))
\end{eqnarray*}
Using Theorem \ref{invariance2} on the torsion of filtered maps we see
that:
\begin{eqnarray*}
\Phi([L_1],[L_2]) &=& \newt(C(\widetilde{L_1}) \to
C(\widetilde{L_2})\\
&=& \newt(\bG_*(C(\widetilde{L_1}))\to \bG_*(C(\widetilde{L_2}))) \\
&=&p^!\newt(C(\widetilde{K_1}) \to C(\widetilde{K_2}))
\end{eqnarray*}
In the general case we choose a common subdivision $(L_3 \to K_3,\cE)$
with fibre structure $F$.  Then:
\begin{eqnarray*}
\Phi([L_1],[L_2]) &=& \Phi([L_1],[L_3]) + \Phi([L_3],[L_2]) \\
&=&p^!\Phi([K_1],[K_3]) + p^!\Phi([K_3],[K_2])\\
&=&p^!\Phi([K_1],[K_2])
\end{eqnarray*}
\end{proof}

\begin{proof}[The proof of Proposition \ref{propone}]
We choose a torsion structure $K\to B$ representing $[B]$.  Then by
Theorem \ref{realization} we may find a fibre bundle torsion structure
$(L \to K,\cE)$ with the given fibre class.  Let $(L' \to K',\cE')$ be
another choice of fibre
bundle torsion structure.  Then by Proposition
\ref{bundle_structure_base_change} we have
\[
\Phi([L],[L']) = p^!([K],[K']) = p^!([B],[B]) = 0
\]
Therefore $[L]$ and $[L']$ are equivalent.
\end{proof}

\begin{proof}[The proof of Proposition \ref{proptwo}]
Choose representatives $(L_1 \to K_1,\cE_1)$, $(L_2 \to K_2,\cE_2)$
and $(L_3 \to K_2,\cE_3)$ for $\fS_{[B_1][F_1]}$,$\fS_{[B_2][F_1]}$
and $\fS_{[B_2][F_2]}$ respectively.  Then
\begin{eqnarray*}
\Phi(\fS_{[B_1][F_1]},\fS_{[B_2][F_2]}) &=& \Phi((L_1 \to
K_1,\cE_1),(L_2 \to K_2,\cE_2)) \\
&&+ \Phi((L_2 \to
K_2,\cE_2),(L_3 \to K_2,\cE_3)) \\
&=& p^!\Phi([B_1],[B_2]) +
  \chi(B)\Phi([F_1],[F_2])
\end{eqnarray*}
by Corollaries \ref{cor:fibre} and \ref{cor:base}.
\end{proof}

\begin{proof}[The proof of Proposition \ref{propthree}]
Let $(L \to K,\cE)$ be a PFBTS representing
  $\fS_{[B][F]}$, and recall that its dual torsion structure $\cD\fS_{[B][F]}$
  is represented by $[\cD L]$, where $\cD \varphi^L_0\colon C(\widetilde{\cD L}) \cong C(\widetilde{L})^{n+k-*}$, as signed, based complexes (see Definition \ref{definition:dual_structure}). 
 We will first  construct a PFBTS representing 
$\fS_{[\cD B][\cD F]}$, and then compare these structures.

 Note first that  $\cD L \to \cD K$ is a cellular map and we
  have a dual PTS for $C(\cD K)$.
We choose a geometric basis on each subcomplex $p^{-1}(\cD \beta) \subset \cD L$ such that
 the composition
 \eqncount
\begin{equation}
\label{eqn:dual}
C(\widetilde{\cD L}) \xrightarrow{(\cD \varphi^L_0)^{-1}} C(\widetilde{L})^{n+k-*}
  \xrightarrow{\theta_{C(\widetilde{L})}} F_*^{dual}C(\widetilde{L})
\end{equation}
is the based identity map of signed \emph{filtered} complexes. 

We now
have an identification
\[\bG_*(C(\widetilde{\cD L})) \to \bG_*(C(\widetilde{L}))^{k-*}\]
 by taking the associated
complex of the above composition, and applying Lemma \ref{kdual} to
identify the 
associated complex of the filtered dual $F^{dual}_*C(\widetilde{L})$
with the $k$-dual $\bG_*(C(\widetilde{L}))^{k-*}$ of the associated complex.
We have a chain equivalence
\eqncount
\begin{equation}
\label{eqn:dual_fbts}
(\cE^{k-*})^{-1}\colon \bG_*(C(\widetilde{\cD L})) \to 
(J^{C(\widetilde{F})}C(\widetilde{K}))^{k-*} =
J^{C(\widetilde{F}))^{n-*}}C(\widetilde{K})^{k-*}
\end{equation}
and this gives an identification
\[(\cE^{b_0})^{k-*}\colon  C(\widetilde{F})^{n-*} \to C(\cD p^{-1}(\cD b_0))\]
However, we may further identify $C(\cD p^{-1}(\cD b_0))$ with $C(\cD F) = C(\cD
p^{-1}(b_0))$ using Remark \ref{remark:transport}, since $b_0$ is the
barycentre of $\cD b_0$.  Then the resulting chain equivalence
\[(\cE^{b_0})^{k-*} = \cD \varphi^F_0
\colon  C(\widetilde{F})^{n-*} \to C(\widetilde{\cD F})\]
coincides with the original identification of  $C(\cD F)$ with
$C(\widetilde{F})^{n-*}$, given by matching cells to dual cells on $F$. This uses our assumption that
 the orientations of $F$, $E$ and $B$ are compatible.

The chain equivalence in (\ref{eqn:dual_fbts}) doesn't yet give
 a PFBTS with underlying
triangulation $\cD L$, since the dual cell complex doesn't have a
base-point.  Let $p'\colon L' \to K'$
be a sub-division of $\cD L \to \cD K$ formed by taking the
barycentric subdivision over the cell $\cD b_0$ only.  In this new
complex we have a base-point $b_0$ in $K'$.
Write $F'$ for  $(p')^{-1}(b_0)$.  

We choose a
PTS on $C(\widetilde{F'})$ such that the subdivision chain equivalence
\[ C(\cD p^{-1}(b_0)) \to C(\widetilde{F'})\]
has trivial absolute torsion, so that $[F'] = [\cD F]$. Similarly, we choose a PTS on $C(\widetilde{K'})$ such that the subdivision
chain equivalence
\[C(\widetilde{\cD K}) \to C(\widetilde{K'})\]
has trivial absolute torsion, so $[K'] = [\cD K]$. The subdivision equivalence 
 is the based identity map away from $\cD b_0$.  
 We now construct a PFBTS $(L' \to K', \cE')$ such that the diagram:
 \eqncount
\begin{equation}\label{diagg}
\vcenter{\xymatrix{\bG_*(C(\widetilde{\cD L})) \ar[d]_{(\cE^{k-*})^{-1}}
   \ar[rr] & & \bG_*(C(\widetilde{L'})) \ar[d]^{\cE'} \\
J^{C(\widetilde{F})^{n-*}}C(\widetilde{K})^{k-*} \ar[r] &
J^{C(\widetilde{F})^{n-*}}C(\widetilde{K'}) \ar[r] &
J^{C(\widetilde{F'})}C(\widetilde{K'})
}}
\end{equation}
commutes, where the top map is associated to the subdivision equivalence
$\cD L \to L'$.  For cells $\beta \in K'$ which correspond to cells in $\cD
L$ (i.e. those away from $\cD b_0$) we define the map ${\cE'}^\beta$ to
be the composition 
\[ C(\widetilde{F'_\beta})\equiv C(\widetilde{\cD F_\beta}) \xrightarrow{((\cE^{\beta})^{k-*})^{-1}}
C(\widetilde{F})^{n-*} \xrightarrow{\cD \varphi^F_0} C(\widetilde{\cD F}) \equiv C(\widetilde{F'})\]
Once we have done this, there are unique maps ${\cE'}^\beta$ for the
remaining $\beta$ such that the sum $\cE'$ is a chain map (existence
is a modified form of Theorem \ref{realization}, uniqueness follows by
considering boundary pairs).  Moreover, the map ${\cE'}^{b_0}\colon
C(\widetilde{F'}) \to C(\widetilde{F'})$ is the
identity map, so after appropriate choices of geometric bases for
the subcomplexes $F'_\beta$,
we have
constructed a PFBTS on $(L' \to K', \cE')$ with fibre class $[\cD F]$
and base class $[\cD B]$; in other words, this structure represents $\fS_{[\cD
    B][\cD F]}$.
It remains to show that
\[\Phi(\cD \fS_{[B][F]},\fS_{[\cD B][\cD F]}) = 0\]
which is equivalent to showing that
\[\Phi(\cD(L \to K,\cE),(L' \to K',\cE')) =0\ .\]

The left, right and bottom maps in diagram (\ref{diagg}) all have trivial absolute torsion
after applying  $i_*\colon K_1^{iso}(\bD\bA(\bZ[\pi_1E])) \to
K_1(\bZ[\pi_1E])$, therefore so does the top map. On the other hand, 
by (\ref{eqn:dual_fbts}) and Proposition \ref{filtered_dual} we have
$$\newt( \cD \varphi^L_0\colon C(\widetilde{L})^{n+k-*} 
\to C(\widetilde{\cD L})) = 
  \newt( \theta_{C(\widetilde{L})}) = 0\ .$$
It follows that
  \[
\begin{array}{l}
\Phi(\cD(L \to K,\cE),(L' \to
   K',\cE')) \\[1ex]
\hspace{1cm}\begin{array}{l}
=\newt(C(\widetilde{L'}) \to
C(\widetilde{L})^{n+k-*})\\[1ex]
= \newt(C(\widetilde{L'}) \to C(\widetilde{\cD L})) +
   \newt(C(\widetilde{\cD L}) \xrightarrow{(\cD \varphi^L_0)^{-1}} C(\widetilde{L})^{n+k-*}) \\[1ex]
=i_*\newt(\bG_*(C(\widetilde{L'}))\to
  \bG_*(C(\widetilde{\cD L}))) \\[1ex] 
=0
\end{array}
\end{array}
\]
as required.
\end{proof}

\begin{proof}[The proof of Theorem B]
We can now complete the proof of Theorem B.
Let $[B]$ and $[F]$ be torsion
structures on $B$ and $F$ respectively. We may assume that
$[F]\in \cS_\bullet(F\xrightarrow{q} E)$ is induced by a PTS on the universal covering of $F$, and therefore
$\Phi([\cD F], [F]) = q_*(\newt(F))$ by Theorem \ref{manifold}. Then
\begin{eqnarray*}
\newt(E) &=& \Phi(\cD \fS_{[B][F]},\fS_{[B][F]})\\
&=& \Phi(\fS_{[\cD B][\cD F]},\fS_{[B][F]}) \\
&=&
\Phi(\fS_{[\cD B][F]},\fS_{[B][F]}) + \Phi(\fS_{[\cD B][\cD F]},\fS_{[\cD B][F]})\\
&=&   p^!\Phi([\cD B][B]) + \chi(B)\Phi([\cD F][F])\\
&=& p^!(\newt(B)) + \chi(B)q_*(\newt(F)) \\
&&\hspace{1cm}\in\widehat{H}^{n+k}(\bZ/2;K_1(\bZ[\pi_1E]))
\end{eqnarray*}
as required.
\end{proof}

\section{The Proof of Theorem A\label{seven}}

Our main result on the multiplicativity of signatures will be deduced from
Theorem \ref{signmod4} and the following application of Theorem B,
in which we  compute the reduced absolute torsion $E$. If $\dim B$ is odd, the vanishing of $\sign(E)$ follows immediately
(as remarked in \cite{atiyah1}), so we may assume that $\dim B$ is even.
\begin{theorem}\label{tauE}
 Let $F \xrightarrow{q} E \xrightarrow{p} B$ be a $PL$
bundle of compatibly oriented closed manifolds. If $\dim B$ is even, then the reduced
absolute torsion
$$\newtbar(E) =\chi(F)\cdot\newtbar(B)) + \chi(B)\cdot\newtbar(F))
\ .$$
\end{theorem}
\begin{proof}
The augmentation map induces a direct sum splitting
$$K_1(\bZ[\pi_1(E,e_0)]) \cong  K_1(\bZ)\oplus
 {\widetilde K}_1(\bZ[\pi_1(E,e_0)])$$
 and $\epsilon_\ast$ restricted to the subgroup $K_1(\bZ)$
 is the identity map. We will apply the augmentation map
 $\epsilon_*$ to the formula for $\newt(E)$ given in Theorem B.
 Note that $\epsilon_*\circ q_* = \epsilon_*$ so the second
 term
 $$\chi(B)\cdot\epsilon_*(q_*\newt(F)) = \chi(B)\cdot\newtbar(F))$$
 as required.

 To evaluate the first term, we need to compute
 $$\epsilon_*(p^!(\newt(B))) = \epsilon_*(p_\ast p^!(\newt(B)))\ .$$
 By Proposition
 \ref{manifold-torsion} we know that $\newt(B) = \tau(\pm h^2)$
 for some $h \in \pi_1(B,b_0)$.
 But the map
 $\phi\colon \pi_1(B,b_0) \to K_1(\bZ)$ defined by
 $\phi(g) = \epsilon_*(p^!(\tau(g)))$
 is a group homomorphism.
 Since $K_1(\bZ) =\cy 2$,  $\phi$ vanishes on
 squares of group elements in $\pi_1(B,b_0)$.  Therefore
 $$\epsilon_\ast(p^!(\newt(B))) = \epsilon_\ast(p^!(\tau(\pm 1)))\ .$$
 However, $\tau(\pm 1)$ lies in the image of the map  $i_\ast\colon K_1(\bZ) \to K_1(\bZ[\pi_1(B,b_0)])$ induced by inclusion.
 By pulling back the bundle
 $E\xrightarrow{p} B$ over the base point $b_0\in B$, and applying
 the naturality formula for pullbacks (see \cite[Cor.~5.3]{lueck1}),
 we obtain the relation
 $$p^!(\newt(B)) =
 i_\ast(p_0^\ast(\tau(\pm 1)))=i_\ast(p_0^\ast(\tau(\pm 1)))$$
 by comparison with the trivial bundle
 $F \times \{b_0\} \xrightarrow{p_0} \{b_0\}$.
 But for a trivial bundle  $p_\ast\circ p^!$ is just multiplication
 by $\chi(F)$, by an easy special case of \cite[Thm.~7.1]{lueck2}.
 Therefore
 $$\epsilon_\ast(p^!(\newt(B))) = \epsilon_\ast(p_\ast(p^!(\newt(B))))=
 = \chi(F)\cdot
 \newtbar(B)$$
 as required, since $\epsilon_*\circ i_\ast$ is the identity map.
\end{proof}
\begin{corollary} If $\dim B$ is even, then
$\newtbar(E) = \newtbar(F\times B)$.
\end{corollary}
\begin{proof}
We apply the product formula, Proposition \ref{product}, and  the
formula just proved.
\end{proof}
This result says that
 the reduced absolute torsion of $E$ is the same as
that of the total space $F\times B$ for the trivial bundle,
provided $\dim B$ is even,
We can now conclude
that the signature $\sign (E)$   agrees with $\sign(F\times B) =
\sign(F)\cdot\sign(B)$ modulo $4$.
\begin{enumerate}
\item If $\dim B = 4j$, $\dim F = 4l$, and $\dim E = 4(j+l) = 4k$
$$\begin{array}{ll}
  \hfill \sign(E) &= 2\newtbar(E) +(2k+1)\chi(E)\\[1ex]
  &=2(\newtbar(B)\chi(F) + \chi(B)\newtbar(F)) + (2k+1)\chi(B)\chi(F)\\[1ex]
   \hfill\sign(B)\sign(F) &= \big [ 2\newtbar(b) +(2j+1)\chi(B)\big ]\big [2\newtbar(F) +(2l+1)\chi(F)\big ]
    \end{array}$$
   and these agree modulo 4.\\
\item If $\dim B = 4j+2$, $\dim F = 4l-2$, and $\dim E = 4(j+l) = 4k$
$$\sign(E)
  =2(\newtbar(B)\chi(F) + \chi(B)\newtbar(F)) + (2k+1)\chi(B)\chi(F)$$ and $\chi(F) \equiv \chi(B) \equiv 0 \Mod 2$, so $\sign(E)\equiv 0 \Mod 4$.\\
  \item If $\dim B$ is even, but $\dim F$ is odd, then both sides are zero.

\end{enumerate}

\noindent
This completes the proof of Theorem A.
\hfill\qed

\section{Fibrations of PD spaces\label{eight}}
The conjecture of \cite{kt1} is stated in a more general situation:
for fibrations of Poincar\'e duality spaces. We don't know yet if the
signature is multiplicative mod 4 for such fibrations, and hope to return to this topic in a future paper.
We will just check that there are
no counter-examples arising from finite coverings.
A formula for the Whitehead torsions in a fibration of Poincar\'e duality spaces was given in \cite{pedersen1}.

Recall that C. T. C. Wall constructed
examples of finite coverings of oriented Poincar\'e complexes
$X$, with $\pi_1(X) = \cy p$, $p$ prime,  and the property that $\sign(\widetilde X) \neq p\cdot \sign (X)$. On the other hand,
the algebraic theory of surgery gives the following congruence.
\begin{lemma}\label{mod8}
 Let $(f,b)\colon Y \to X$ be a degree one normal
map of oriented, finite Poincar\'e duality spaces. Then
$\sign(Y) - \sign(X) \equiv 0 \Mod 8$.
\end{lemma}
\begin{proof}
By \cite[p.~229]{ra11}
any such degree one normal map  has a quadratic signature
$\sigma_*(f,b)$, with the property that
$$\sign(\sigma_*(f,b))= \sign(\sigma^*(Y)) - \sign(\sigma^*(X)) $$
and $\sign(X) = \sign(\sigma^*(X))$ for any Poincar\'e complex.
On the other hand, $\sign(\sigma_*(f,b)) \equiv 0 \Mod 8$ since
the signature of an even, unimodular symmetric bilinear form over
$\bZ$ is divisible by 8 (see \cite[3.11]{hnk1}).
\end{proof}
\begin{corollary} Let $X'\to X$ be a finite covering of degree $d$
of oriented, finite Poincar\'e complexes. Then $\sign(X') - d\cdot\sign(X)
\equiv 0 \Mod 8$.
\end{corollary}
\begin{proof}
For any finite covering
$X' \to X$ of oriented, finite Poincar\'e duality spaces, we can construct a degree
one normal map
$$(f,b)\colon X' \sqcup -(d-1) X \to X\ .$$
Then Lemma \ref{mod8} implies that
$$\sign(X') - d\cdot\sign(X) = \sign(X' \sqcup -(d-1) X) - \sign(X) \equiv 0 \Mod 8\ .$$
\end{proof}

\section{Filtered chain complexes\label{three}}

In this section, we give a self-contained  treatment of the
absolute torsion of signed filtered chain complexes. These results
are the algebraic foundation for the torsion calculations in this
paper, but this section  can be read independently of the previous
sections.
 The main results are the Invariance Theorems \ref{invariance1}, \ref{invariance2}, generalizing results previously obtained for the reduced torsion by Milnor \cite{milnor1}, Maumary \cite{maumary}, and Munkholm \cite{munkholm}.  These results
 express the absolute torsion $\newt(f) \in K_1(\bA)$
of a filtered chain equivalence $f\colon C \to D$ of signed filtered
complexes in terms of the filtration quotients.

\subsection{Filtered complexes}\label{filter}
We first need some notation and definitions.
\begin{definition} \label{fil}
Let $\bA$ be an additive category.
\begin{enumerate}
\item A \emph{$k$-filtered object} $F_*M$ in $\bA$ is an object $M$
in $\bA$ together with a direct sum decomposition
$$M=M_0\oplus M_1 \oplus \dots \oplus M_k$$
which we regard as a length $k$ filtration
$$F_{-1}M=0 \subseteq F_0M \subseteq F_1M \subseteq \dots
\subseteq F_kM=M$$
with
$$F_jM=M_0 \oplus M_1 \oplus \dots \oplus M_j~~(0 \leqslant j \leqslant k)~.$$
\item  A \emph{filtered morphism} $f\colon F_*M \to F_*N$ of $k$-filtered
objects in $\bA$ is a morphism in $\bA$ of the type
$$\begin{array}{l}
f=\begin{pmatrix} f_0 & f_1 & f_2 & \dots & f_k \\
0 & f_0 & f_1 & \dots & f_{k-1} \\
0 & 0 & f_0 & \dots & f_{k-2} \\
\vdots & \vdots & \vdots & \ddots & \vdots \\
0 & 0 & 0 & \dots & f_0 \end{pmatrix}~:~ M=\bigoplus\limits_{s=0}^k M_s \to
N=\bigoplus\limits_{s=0}^k N_s
\end{array}$$
so that
$$f(F_jM) \subseteq F_jN~~(0 \leqslant j \leqslant k)~.$$
The  $(u,v)$-component  of this upper triangular matrix is a
morphism $f_{v-u}\colon M_{v} \to N_{u}$, $0\leqslant u\leqslant
v\leqslant k$, where
$f_j\colon M_* \to N_{*-j}$, $0\leqslant j \leqslant k$,   are graded morphisms in $\bA$. \\
\item  A \emph{$k$-filtered complex} $F_*C$ in $\bA$ is a finite chain complex
$C$ in $\bA$ with $k$-filtered objects $C_r=\bigoplus_{s=0}^k C_{r,s}$ such that the differentials
$d\colon F_*C_r \to F_*C_{r-1}$ are filtered morphisms. The matrix components of $d$ are  maps $d_j\colon C_{r,s} \to C_{r-1,s-j}$. \\
\item A \emph{filtered chain map} $f\colon F_*C \to F_*D$ is a chain map
$f\colon C \to D$ such that  $f_r\colon F_*C_r \to F_*D_r$ is a filtered morphism in each  degree. The component maps have the form
$f_j\colon C_{r,s} \to D_{r,s-j}$.\\
\item A \emph{filtered chain homotopy} $g\colon f\simeq f'\colon F_*C \to F_*D$
between filtered chain maps $f,f'$ is a collection
$\{g\colon C_r \to D_{r+1}\,\vert\,r \in \bZ\}$ of morphisms in $\bA$ such that
$$f-f'=dg+gd\colon C_r \to D_r$$
and
$$\begin{array}{l}
g=\begin{pmatrix} g_0 & g_1 & g_2 & \dots & g_k \\
g_{-1} & g_0 & g_1 & \dots & g_{k-1} \\
0 & g_{-1} & g_0 & \dots & g_{k-2} \\
\vdots & \vdots & \vdots & \ddots & \vdots \\
0 & 0 & 0 & \dots & g_0 \end{pmatrix}
\end{array}$$
has component maps
$$g_j\colon C_{r,s}\to  D_{r+1,s-j}~~(-1 \leqslant j \leqslant k)~.$$
\item A \emph{filtered contraction} of a $k$-filtered complex $F_*C$ in $\bA$
is a filtered chain homotopy
$$\Gamma\colon 1~\simeq~0\colon F_*C \to F_*C~.$$
\item The \emph{filtered mapping cone} of a filtered chain map $f\colon C \to D$
of $k$-filtered complexes is the $(k+1)$-filtered complex
$F_*\FC(f)$ with
$$\FC(f)_{r,s}=D_{r,s} \oplus C_{r-1,s-1}~~(0 \leqslant s \leqslant k+1)$$
and the differential $d^{\,\FC(f)}\colon \FC(f)_r \to
\FC(f)_{r-1}$ is in upper triangular block form with $2\times 2$
block entries ($0\leqslant j \leqslant k+1$):
$$d^{\,\FC(f)}_j = \left (\vcenter
{\xymatrix@C-2pc@R-2pc{d^D_j&(-)^{r-1}f_{j-1}\\0&d^C_j} }
\right )\colon D_{r,s} \oplus C_{r-1,s-1} \to D_{r-1,s-j} \oplus C_{r-2,s-j-1}$$
\end{enumerate}
\hfill\qed
\end{definition}
\noindent
The rearrangement map
$$\FC(f)_r= \sum_{s=0}^k \left (D_{r,s}\oplus C_{r-1,s-1}\right ) \xrightarrow{\cong}
\sum_{s=0}^k D_{r,s} \oplus \sum_{s=1}^k  C_{r-1,s-1} = \cC(f)_r$$
define an isomorphism of unfiltered chain complexes $
 \rho\colon \FC(f)\xrightarrow{\cong} \cC(f)$.

\begin{theorem}\label{contract}
A filtered chain map $f\colon F_*C \to F_*D$ of $k$-filtered
chain complexes is a filtered chain equivalence if and only if the
filtered mapping cone $F_*\FC(f)$ is filtered contractible.
\end{theorem}
\begin{proof} By definition, $f$ is a filtered chain equivalence
if and only if there exist a filtered chain map $g\colon F_*D \to F_*C$ and
filtered chain homotopies
$$\Gamma\colon 1~\simeq~gf\colon F_*C \to F_*C~,~
\Delta\colon 1~\simeq~fg\colon F_*D \to F_*D~.$$
\indent A filtered contraction $e\colon 1 \simeq 0\colon F_*\FC(f) \to F_*\FC(f)$
is of the type
$$e=\begin{pmatrix} \Delta & h \\ (-)^rg & \Gamma \end{pmatrix}\colon
\FC(f)_r~\cong~D_r \oplus C_{r-1} \to
\FC(f)_{r+1}~\cong~D_{r+1} \oplus C_r$$
with $g\colon F_*D \to F_*C$ a filtered chain map and
$\Gamma\colon 1 \simeq gf$, $\Delta\colon 1\simeq fg$ filtered chain homotopies,
so that $f$ is a filtered chain equivalence.\\
\indent Conversely, suppose that $f$ is a filtered chain equivalence.
The morphisms defined by
$$e'=\begin{pmatrix} \Delta & 0 \\ (-)^rg & \Gamma \end{pmatrix}\colon
\FC(f)_r~\cong~D_r \oplus C_{r-1} \to
\FC(f)_{r+1}~\cong~D_{r+1} \oplus C_r$$
are such that the morphisms
$$\begin{array}{l}
de'+e'd=\begin{pmatrix} 1 & (-)^r(f\Gamma-\Delta f)\\
0 & 1 \end{pmatrix}~\colon \\[2ex]
\hskip100pt
\FC(f)_r~\cong~D_r \oplus C_{r-1} \to
\FC(f)_r~\cong~D_r \oplus C_{r-1}
\end{array}$$
are automorphisms. The morphisms
$$e=(de'+e'd)^{-1}e'\colon \FC(f)_r~\cong~D_r \oplus C_{r-1} \to
\FC(f)_{r+1}~\cong~D_{r+1} \oplus C_r$$
define a filtered contraction $e\colon 1 \simeq 0\colon F_*\FC(f) \to F_*\FC(f)$.
\end{proof}

\subsection{The associated complex}
We now describe how to pass from a filtered complex to its associated graded complex. The differentials $d\colon C_r \to C_{r-1}$
of a $k$-filtered complex $F_*C$ satisfy  the matrix identities
$$d^2=0\colon C_r \to C_{r-2}\ .$$  
By carrying out the multiplication $d^2$ as a product of upper triangular matrices, we obtain
a number of recursion formulas among the component maps
$d_j\colon C_{r,s} \to C_{r-1,s-j}$. 
In particular, we have the formulas
$$\begin{array}{l}
(d_0)^2=0\colon C_{r,s} \to C_{r-2,s}~,\\[.6ex]
d_0d_1+d_1d_0=0\colon C_{r,s} \to C_{r-2,s-1}~,\\[.6ex]
(d_1)^2+d_0d_2+d_2d_0=0\colon C_{r,s} \to C_{r-2,s-2}~.
\end{array}$$

\begin{definition}
The \emph{associated complex} of a $k$-filtered complex
$F_*C$ in $\bA$ is the $k$-dimensional chain complex in
the derived category $\bD(\bA)$
$$\bG_*(C) :\ \xymatrix{G_k(C) \to\dots\to G_{r+1}(C) \ar[r]^-{d_*}&
G_r(C) \ar[r]^-{d_*}&  \dots\to G_0(C)}$$
with differential
$$d_*=(-)^{s}d_1\colon G_r(C)_s=C_{r+s,r} \to G_{r-1}(C)_s=C_{r+s-1,r-1}\ .$$
The individual terms $G_r(C)$ are objects in $\bD(\bA)$, with `internal' differential
$$d_{G_r(C)}=d_0\colon G_r(C)_s=C_{r+s,r} \to
G_r(C)_{s-1}=C_{r+s-1,r}\ .$$
\hfill\qed
\end{definition}
\begin{example}[Tensor Products]\label{tp_one}
Let $C$ and $D$ be chain complexes in $\bA(R)$ and $\bA(S)$
respectively, where $R$ and $S$ are rings. If $\dim C =k$, then
the tensor product complex $C\otimes D$ (over $\bZ$) admits a
$k$-filtered  structure with $(C\otimes D)_{r,s} = C_s\otimes
D_{r-s}$ for $0\leqslant s \leqslant k$. The filtered differential
is defined by
\begin{eqnarray*} d_0&=& 1\otimes d^D\colon C_s\otimes D_{r-s}
\to C_s\otimes D_{r-s-1}\\
 d_1&=&(-)^{r-s}d^C\otimes 1\colon C_{s-1}\otimes D_{r-s}
\to C_{s-1}\otimes D_{r-s}
\end{eqnarray*}
 and $d_j = 0$ for $j\geqslant 2$. The associated complex has $G_r(C\otimes D)_s = C_r\otimes D_s$
and differential $d_* = d^C\otimes 1\colon C_r\otimes D_s
\to C_{r-1}\otimes D_s$.
\hfill\qed
\end{example}
When we work with the associated complex, it is useful to translate the maps $d_j\colon C_{r,s} \to C_{r-1,s-j}$ into the new notation, so that $$d_j\colon G_r(C)_s \to G_{r-j}(C)_{s+j-1}~.$$
The relation $d_0d_1+d_1d_0=0$ then implies that $d_*\colon G_r(C) \to G_{r-1}(C)$ is a chain map. The relation $(d_1)^2+d_0d_2+d_2d_0=0$ shows that
$(d_1)^2$ is chain homotopic to zero, and hence $(d_*)^2=0$ in the derived category.
\begin{proposition}\label{filprop}
\hfill
\begin{enumerate}
\item A filtered chain map $f\colon F_*C \to F_*D$ induces a chain map
in $\bD(\bA)$
$$\bG_*(f)\colon \bG_*(C) \to \bG_*(D)~.$$
\item A filtered chain homotopy $g\colon f \simeq f'\colon F_*C \to F_*D$
induces a chain homotopy in $\bD(\bA)$
$$\bG_*(g)\colon \bG_*(f) \simeq \bG_*(f')\colon \bG_*(C) \to \bG_*(D)~.$$
\item A filtered chain equivalence $f\colon F_*C \to F_*D$ induces
a chain equivalence $\bG_*(f)\colon \bG_*(C) \to \bG_*(D)$ in $\bD(\bA)$.\\
\item A filtered chain contraction $\Gamma\colon 1 \simeq 0\colon F_*C \to F_*D$ induces
a chain equivalence $\bG_*(\Gamma)\colon 1 \simeq 0\colon \bG_*(C) \to \bG_*(D)$ in $\bD(\bA)$.
\end{enumerate}
\end{proposition}
\begin{proof}
(i) For a filtered chain map $f\colon F_*C \to F_*D$ the identities
$$d f =f d\colon C_r \to D_r$$
expressed in upper triangular matrix form lead to relations involving the component maps $f_i$ and $d_j$.
In particular we obtain
$d_0f_0=f_0d_0\colon C_{r,s} \to D_{r-1,s}$
 and the map
$G_r(f)=f_0\colon G_r(C) \to G_r(D)$ is a chain map. The relation
$d_1f_0 -f_0d_1 = f_1d_0-d_0f_1$ shows that $(-)^{s}f_1\colon
G_r(C)_s \to G_{r-1}(D)_{s+1}$ gives a chain homotopy
$d_*f_0 - f_0 d_* \simeq 0$, and so $f_0$ gives a chain map
$\bG_*(f)\colon\bG_*(C) \to \bG_*(D)$ in $\bD(\bA)$.

(ii) For a filtered chain homotopy $g\colon f\simeq f'\colon F_*C \to F_*D$ the identities
$$f-f'=dg+gd\colon C_r \to D_r$$
expressed  in upper triangular matrix form again lead to various relations.
In particular
$$\begin{array}{l}
d_0g_{-1}+g_{-1}d_0=0\colon C_{r,s} \to D_{r,s+1}~,\\[.6ex]
f_0-f'_0=d_0g_0+g_0d_0~+d_1g_{-1}+g_{-1}d_1\colon C_{r,s} \to D_{r,s}~.
\end{array}$$
In the $\bG_*$-notation,  we have maps $$g_j\colon G_r(C)_s \to G_{r-j}(D)_{s+j+1}~.$$
Since $g_0\colon G_r(C)_s \to G_{r}(D)_{s+1}$ gives a chain null-homotopy $d_0g_0+g_0d_0\simeq 0\colon G_r(C) \to G_{r}(D)$, the map
$$\bG_*(g)=(-)^{s}g_{-1}\colon G_r(C)_s \to G_{r+1}(D)_s$$
gives a null-homotopy $f_0 - f'_0\simeq 0$ in the derived category
$\bD(\bA)$.
Parts (iii) and (iv)  follow immediately  from parts (i) and (ii).
\end{proof}

\begin{theorem} \label{equivalence}
{\rm (i)} A $k$-filtered complex $F_*C$ in $\bA$ is filtered contractible if and only
if the associated complex $\bG_*(C)$ is contractible in $\bD(\bA)$.\\
{\rm (ii)} A filtered chain map $f\colon F_*C \to F_*D$ of $k$-filtered
chain complexes is a filtered chain equivalence if and only if the
associated chain map $\bG_*(f)\colon \bG_*(C) \to \bG_*(D)$ is a chain equivalence in $\bD(\bA)$.
\end{theorem}
\begin{proof} (i) A filtered contraction $\Gamma\colon 1\simeq 0\colon F_*C \to F_*C$
induces a contraction $\bG_*(\Gamma)\colon 1 \simeq 0\colon \bG_*(C) \to \bG_*(C)$
by Proposition \ref{filprop} (iv).\\
\indent Conversely, suppose given a contraction of $\bG_*(C)$ in $\bD(\bA)$
$$e\colon 1~\simeq~0\colon \bG_*(C) \to \bG_*(C)~,$$
as represented by chain maps $e\colon (G_r(C),d_0) \to (G_{r+1}(C),d_0)$,  for which
there exist chain homotopies
$$h\colon 1  ~\simeq~d_*e+ed_*\colon G_r(C) \to G_r(C)$$
with
$$(-)^{s}(d_1e+ed_1)+d_0h+hd_0=1\colon
G_r(C)_s=C_{r+s,r} \to G_r(C)_s=C_{r+s,r}~.$$
In terms of the filtration on $C$, the maps $h\colon G_r(C)_s \to G_r(C)_{s+1}$ and $e\colon G_r(C)_s \to G_{r+1}(C)_s $ give maps
$h_s\colon C_{r,s} \to C_{r+1,s}$ and $e_{r,s}\colon C_{r,s} \to C_{r+1,s+1}$.
The relation given by $h\colon 1\simeq d_* e + ed_*$ gives
$$(-)^{r+s}(d_1 e + e d_1) + d_0 h + h d_0=1\colon C_{r,s} \to C_{r,s}$$
Let  $\hat e_{r,s}=(-)^{r-1}e$, and define morphisms
$$\begin{array}{l}
\beta=\begin{pmatrix}
h & 0 & 0 & \dots & 0\\
\hat e & h & 0 & \dots & 0 \\
0 & \hat e & h & \dots & 0 \\
\vdots & \vdots & \vdots & \ddots \\
0 & 0 & 0 & \dots & h
\end{pmatrix}
\colon  C_r \to C_{r+1}
\end{array}$$
The morphisms $\alpha\colon C_r \to C_r$ defined by $$\alpha=d\beta+\beta d\colon C_r \to C_r$$ are in upper triangular form, with the identity  on the main diagonal.
For the components in the $(s+1,s)$ positions, note that
we get $(-)^{r-1}d_0 e + (-)^r e d_0 =0$ since $e$ a chain map. In the diagonal $(s,s)$ positions we get the sum of composites
$$\xymatrix@R-15pt{
C_{r,s} \xrightarrow{(-)^{s} d_1} C_{r-1,s-1} \xrightarrow{(-)^{r} e} C_{r,s}&
C_{r,s} \xrightarrow{(-)^{r-1} e} C_{r+1,s+1} \xrightarrow{(-)^{s+1} d_1} C_{r,s}\\
C_{r,s} \xrightarrow{d_0} C_{r-1,s} \xrightarrow{ h} C_{r,s}&
C_{r,s} \xrightarrow{h} C_{r+1,s} \xrightarrow{d_0} C_{r,s}}
$$
which equals
$(-)^{r+s}(d_1  e + e d_1) + d_0 h + h d_0 =1$ by the relation above.

Each $\alpha\colon C_r \to C_r$ is an automorphism, with
$$\begin{array}{l}
d \alpha=d \beta d=\alpha d\colon C_r \to C_{r-1}~,\\[.6ex]
\alpha^{-1}d=d\alpha^{-1}\colon C_r \to C_{r-1}~.
\end{array}$$
The morphisms
$$\Gamma=\beta\alpha^{-1}\colon C_r \to C_{r+1}$$
are such that
$$d\Gamma+\Gamma d=(d\beta +\beta d)\alpha^{-1}=1\colon C_r \to C_r~,$$
and define a filtered contraction of $F_*C$
$$\Gamma \colon 1~\simeq~0\colon C \to C~.$$
(ii) By Theorem \ref{contract} $f$ is a filtered chain equivalence
if and only if $F_*\FC(f)$ is filtered contractible. By (i)
$F_*\FC(f)$ is filtered contractible if and only if
$\bG_*(\FC(f))=\cC(G_*(f))$ is chain contractible.
\end{proof}

%%%%%%%%%%%%%%%%%%%%%%%%%%%%%
\subsection{Splitting and Folding}\label{zero_b}
 In an additive category $\bA$ we don't necessarily have kernels and cokernels, but we can define split exact sequences.
A \emph{direct sum system}
$\ A \xymatrix{\ar@<-.6ex>@{<-}[r]_-{\Gamma}
\ar@<.6ex>[r]^-{f} &} B \xymatrix{\ar@<-.6ex>@{<-}[r]_-{\Delta}
\ar@<.6ex>[r]^-{g} &} C \ $
 is a collection of morphisms  in $\bA$ such that
$$(\Delta~f)\colon C \oplus A \to B~,~{g\choose \Gamma}\colon
B \to  C \oplus A$$
are inverse isomorphisms in $\bA$. Then we say that
a morphism $f\colon A \to B$ in $\bA$ is a \emph{split injection} (respectively $g\colon B \to C$ is  a \emph{split surjection}) if the morphism
extends to a direct sum system. Then
$$0 \to A \xymatrix{\ar[r]^-{\di{f}}&} B \xymatrix{\ar[r]^-{\di{g}}&} C \to 0$$
is a \emph{short exact sequence} in $\bA$ if $f$ and $g$ can be extended to a direct sum system.
A short exact sequence of chain complexes in $\bA$ is one which
is short exact in each degree.

There is a useful criterion for the existence of a split injection in the derived category $\bD(\bA)$, expressed in terms of any representative chain map $f$ for a given morphism.
\begin{proposition} \label{derivedsplit} A morphism $f\colon C \to D$ in $\bD(\bA)$ is a split
injection if and only if
there exists a chain map $\Gamma\colon D \to C$ and a chain homotopy
$h\colon \Gamma f \simeq 1\colon C \to C$. If  $(\Gamma, h)$ exists, then there is a direct sum system
$$C \xymatrix{\ar@<-.6ex>@{<-}[r]_-{\di{\Gamma}}
\ar@<.6ex>[r]^-{\di{f}} &} D \xymatrix{\ar@<-.6ex>@{<-}[r]_-{\di{\Delta}}
\ar@<.6ex>[r]^-{\di{g}} &} \cC(f) $$
in $\bD(\bA)$,
with $g={1\choose 0}\colon D_r \to \cC(f)_r=D_r\oplus C_{r-1}$ and
$$\Delta=(1-f\Gamma\ \ (-)^{r+1}fh)\colon\, \cC(f)_r=D_r\oplus C_{r-1} \to D_r\ .$$
\end{proposition}
\begin{proof} Assuming there exist such $\Gamma,h$ we define a direct sum
system in $\bD(\bA)$ by the given formulas.
The chain maps
$$(\Delta~f)\colon \cC(f) \oplus C \to D~,~
{g\choose \Gamma}\colon D \to \cC(f) \oplus C$$
are inverse chain equivalences,
and the morphisms
$$e=\begin{pmatrix} 0 & 0 & 0 \\
(-)^r\Gamma & -h & (-)^{r+1} \\
h\Gamma & (-)^{r+1}h^2 & h \end{pmatrix}\colon D_r \oplus C_{r-1} \oplus C_r
\to D_{r+1} \oplus C_r \oplus C_{r+1}$$
define a chain homotopy
$$e\colon {g\choose \Gamma}(\Delta~f)~\simeq~1\colon
\cC(f) \oplus C \to \cC(f) \oplus C~.$$
The converse is clear.
\end{proof}

\begin{corollary}\label{tech}
Let $f\colon C \to D$ be a split injection in $\SD(\bA)$. Then
$$i_*\newt({g\choose \Gamma}) = i_*\newt(\Delta~ f) = 0 \in K_1^{iso}(\bA)\ .$$
\end{corollary}
\begin{proof}
This follows immediately from the last result and the formula in \cite[Prop. 13.6]{ak1}, since $\newt(f) = \newt(\cC(f))$.
\end{proof}

We will need a variant of the folding construction used by Whitehead
\cite{whitehead}  to define the torsion of a contractible complex.
Let
$$C~\colon ~ C_k \xymatrix{\ar[r]^-{\di{d}}&} C_{k-1}
\xymatrix{\ar[r]^-{\di{d}}&} C_{k-2}  \to \dots \to C_0$$
be a $k$-dimensional chain complex in $\bA$,
such that the boundary map
$d\colon C_k \to C_{k-1}$ extends to a direct sum system
$$C_k \xymatrix{\ar@<-.6ex>@{<-}[r]_-{\di{\Gamma}}
\ar@<.6ex>[r]^-{\di{d}} &} C_{k-1}
\xymatrix{\ar@<-.6ex>@{<-}[r]_-{\di{\Delta}}
\ar@<.6ex>[r]^-{\di{g}} &} C'_{k-1}$$
We use this direct sum system to define the $(k-1)$-dimensional chain complex
$$C':\  \xymatrix{C'_{k-1}\ar[r]^{\di{d\Delta}}& C_{k-2}\ar[r]^{\di{d}}&
C_{k-3}\ar[r]&  \dots \ar[r]& C_0}
$$
in $\bA$, called the {\it abelian $k$-folding} of $C$. An
\emph{elementary} chain complex in $\bA$ is a contractible complex
with non-zero chain groups only in two adjacent degrees.

\begin{proposition} \label{elementary} Let $C'$ be the abelian
$k$-folding of $C$, and let $E$ be the elementary chain complexes
with $C_k \xrightarrow{1} C_k$ in adjacent degrees $(k,k-1)$.
Then there exist an isomorphism
${i\choose j}\colon C ~\cong~ C' \oplus E$
with isomorphism torsion
$$
\tauiso({i\choose j})=
(-)^{k-1}\tauiso({g\choose \Gamma}\colon
C_{k-1} \to C'_{k-1} \oplus C_k) \in K^{iso}_1(\bA)~.
$$
\end{proposition}
\begin{proof}
The chain isomorphism ${i\choose j}\colon C \to C'\oplus E$ is defined by the diagram
$$\xymatrix@C+8pt{
C_k\ar[r]^d\ar@{=}[d]& C_{k-1}\ar[r]^d\ar[d]^{g\choose\Gamma}& C_{k-2}\ar[r]^d\ar@{=}[d]& \dots \ar[r]& C_0\ar@{=}[d]\\
C_k\ar[r]^(.4){0\choose 1}& C'_{k-1}\oplus C_k\ar[r]^(.6){(d\Delta~0)}& C_{k-2}\ar[r]^d& \dots \ar[r]& C_0
}$$
and the relations $gd=0$, $\Gamma d=1$.
\end{proof}
%%%%%%%%%%%%%%%%%%%%%%%%%%%%
\subsection{The Invariance Theorem for filtered contractible complexes}
We will first establish the sign conventions for filtered complexes.
\begin{definition}
 A \emph{$k$-filtered signed complex} $(F_*C,\eta_{F_*C})$ in
$\bA$ is a $k$-filtered complex $C$ together with signs
$$\begin{array}{l}
\eta_{\bG_*(C)}\in
{\rm im}(\epsilon\colon K_0(\SD(\bA)) \otimes K_0(\SD(\bA)) \to K_1^{iso}(\SD(\bA)))~,\\[.6ex]
\eta_{G_j(C)}
\in {\rm im}(\epsilon\colon K_0(\bA) \otimes K_0(\bA) \to K_1^{iso}(\bA))~~
(0 \leqslant j \leqslant k)~.
\end{array}
$$
The sign of $F_*C$ and the sign of $C$ (as an unfiltered complex)
are set to be
$$\begin{array}{ll}
\eta_{F_*C}&=~\eta_C=i_*\eta_{\bG_*(C)}+\eta_{G_0C\oplus S(G_1C \oplus S(G_2C\oplus
\dots\oplus SG_kC)\ldots)}\\[.6ex]
&\in {\rm im}(\epsilon\colon K_0(\bA) \otimes K_0(\bA) \to K_1^{iso}(\bA))~.
\end{array}$$
The \emph{associated chain complex} $(\bG_*(C),\eta_{\bG_*(C)})$
is a $k$-dimensional chain complex in $\SD(\bA)$. We will usually denote this \emph{signed} complex in $\bD(\bA)$ just by $\bG_*(C)$.
\hfill\qed
\end{definition}

The sign conventions have been chosen so that for a filtered contractible
$k$-filtered signed complex $(F_*C,\eta_{F_*C})$ in $\bA$
$$\newt(C,\eta_C)=\newt(F_*C,\eta_{F_*C}) \in K^{iso}_1(\bA)~.$$

The underlying (unfiltered) complex of a filtered complex
has a useful \emph{iterated mapping cone} description.
Given a $k$-filtered chain complex $F_*C$
in $\bA$, we define chain complexes $T_{\ell, r}(C) = S^{-r}(F_\ell C/F_{r-1}C)$ with chain groups
$$T_{\ell, r}(C)_s =  (F_\ell C/F_{r-1}C)_{r+s} = C_{r+s,r} \oplus C_{r+s,r+1} \oplus \dots \oplus C_{r+s,\ell}$$
for $1 \leqslant r \leqslant \ell\leqslant k$. The differential on
$T_{\ell, r}(C) $
is the one induced by $d\colon C \to C$, and the sign
$$\eta_{T_{\ell,r}(C)} =
\eta_{G_rC\oplus S(G_{r+1}C \oplus S(G_{r+2}C\oplus
\dots\oplus SG_\ell C)\ldots)}$$
Note that the definition
$G_r(C)_s = C_{r+s,r}$ gives the formula
$$T_{\ell, r}(C) _s = G_r(C)_s \oplus G_{r+1}(C)_{s-1} \oplus \dots \oplus G_{\ell}(C)_{r+s-\ell}\ .$$
We use this expression  to define  maps
$$\partial_{\ell,r}\colon T_{\ell, r}(C)_s \to G_{r-1}(C)_s\quad
(1 \leqslant j \leqslant k)$$
by the row matrix
$$\partial_{\ell,r}=(-)^s(d_1~d_2~\dots~d_{\ell-r+1})$$
where $d_j\colon G_r(C)_s \to G_{r-j}(C)_{s+j-1}$ are the entries in
the matrix expression for $d$. If $\ell =k$ we will set
$\partial_r:=\partial_{k,r} $ to simplify the notation.
We remark that $G_r(C)=T_{r,r}(C)$
and $\partial_{r,r} = d_*\colon G_r(C) \to G_{r-1}(C)$.

\begin{proposition}\label{mappingcone}
Let $(F_*C, \eta_{F_*C})$ be a $k$-filtered signed complex in $\bA$.
\begin{enumerate}
\item  The maps $\partial_{\ell,r} \colon
T_{\ell, r}(C) \to G_{r-1}(C)$ are chain maps.
\item $T_{\ell,r-1}(C) = \cC(
T_{\ell, r}(C)  \xrightarrow{\partial_{\ell,r}} G_{r-1}(C))$,  for
$1 \leqslant r \leqslant \ell\leqslant k$, as signed complexes.\\
\item There is a short exact sequence
$$ T_{\ell_1,r}(C) \to T_{\ell_2,r}(C) \to S^{\ell_1-r + 1}T_{\ell_2,\ell_1+1}(C)$$
in $\SD(\bA)$,
for $1\leqslant r\leqslant \ell_1 < \ell_2\leqslant k$.\\
\item The signed complex
$(C, \eta_C -i_*\eta_{\bG_*(C)}) = T_{k,0}(C)$
 as an object in $\SD(\bA)$.
 \end{enumerate}
\end{proposition}

\begin{proof} The claim that $\partial_{\ell,r}$ is a chain map follows immediately from the relation $d^2=0$. Part (ii) follows directly from the definition, and part (iii) can be checked inductively.
Note that, according to our conventions,  the sign term
$$\eta_C - i_*\eta_{\bG_*(C)} = \eta_{G_0C\oplus S(G_1C \oplus S(G_2C\oplus
\dots\oplus SG_kC)\ldots)}$$
 is just the sign of the iterated mapping
cone structure on $T_{k,0}(C)$.
\end{proof}

\begin{definition}
Let $(F_*C,\eta_{F_*C})$
be a $k$-filtered signed complex in $\bA$. The  \emph{$(k-1)$-amalgamation}  of $(F_*C,\eta_{F_*C})$ is
 the $(k-1)$-filtered signed complex $(F_*C',\eta_{F_*C'})$ in $\bA$ with filtration summands
$$C'_{r,s}=C_{r,s}, \text{\ if\ } 0\leqslant s<k-1\text{\ and\ } C'_{r,k-1} = C_{r,k-1} \oplus C_{r,k}$$
in each degree $r$.
The differentials $d'_j = d_j$ except in filtration
degree $k-1$, where
$$d'_j = \mmatrix{d_j}{d_{j+1}}{0}{d_j}
\colon C_{r,k-1} \oplus C_{r,k} \to C_{r-1,k-1-j} \oplus C_{r-1,k-j}$$
and
$$\begin{array}{l}
G_r(C')=\begin{cases}
G_r(C)&\text{\ if\ }0 \leqslant r \leqslant k-2 \\
\cC(G_k(C) \xrightarrow{\partial_k} G_{k-1}(C))&\text{\  if\ } r=k-1~.
\end{cases}
\end{array}$$
\smallskip
\noindent
where $\partial_k:=\partial_{k,k} $.
The signs are given by the formulas
$$\begin{array}{l}
\eta_{\bG_*(C')}=\eta_{\bG_*(C)}~,\\[.6ex]
\eta_{G_j(C')}=\begin{cases}
\eta_{G_j(C)}&\hbox{\rm if}~0 \leqslant j \leqslant k-2 \\
\eta_{G_{k-1}(C) \oplus SG_{k}(C)}&\hbox{\rm if}~j=k-1~.
\end{cases}
\end{array}$$
\hfill\qed
\end{definition}

\begin{lemma}\label{fiveten}
$G_{k-1}(C') = \cC(G_k(C)  \xrightarrow{\partial_k}G_{k-1}(C))$ as objects in $\SD(\bA)$.
\end{lemma}
\begin{proof}
The associated $(k-1)$-dimensional chain complex $\bG_*(C')$ in $\SD(\bA)$ has the
differential $d'_*\colon G_r(C') \to G_{r-1}(C')$  given by
$$\begin{array}{l}
d_{\bG_*(C')}=
\begin{cases}
d_{*}\colon G_j(C) \to G_{j-1}(C)&
\text{\ if\ }0 \leqslant j \leqslant k-2 \\
\partial_{k-1}\colon \cC(G_k(C) \xrightarrow{\partial_k} G_{k-1}(C)) \to G_{k-2}(C)&
\text{\ if\ } j=k-1~.
\end{cases}
\end{array}$$
Since $G_{k-1}(C') = G_{k-1}(C)\oplus SG_k(C)$, the signs
on $G_{k-1}(C')$ and $\cC(\partial_k)$ agree.
\end{proof}
Note that $\eta_{G_0C'\oplus S(G_1C' \oplus S(G_2C'\oplus
\dots\oplus SG_{k-1}C')\ldots)}$ is the sign of the new iterated mapping
cone structure on $C$, and so $\eta_{F_*C'}=\eta_{F_*C}=\eta_C$ as
required,

\begin{theorem}[{\bf Invariance}] \label{invariance1}
For a filtered contractible  $k$-filtered signed complex $(F_*C,\eta_{F_*C})$
in $\bA$
$$\newt(C,\eta_C)=i_*\newt(\bG_*(C),\eta_{\bG_*(C)}) \in K^{iso}_1(\bA)~.$$
\end{theorem}
\begin{proof} The proof is by induction on $k$. The result is true for $k=0$,
since in that case $(\bG_*(C),\eta_{\bG*(C)})$ is
concentrated in degree 0, and $\eta_C = i_*\eta_{\bG*(C)} + \eta_{G_0(C)}$. Therefore
$$i_*\newt(\bG_*(C),\eta_{\bG_*(C)}) = \newt(G_0(C),\eta_{G_0(C)}) +
i_*\eta_{\bG*(C)} = \newt(C,\eta_C)\ .$$

So assume that $k \geqslant 1$, and that the
result is true for $(k-1)$. Let $(F_*C',\eta_{F_*C'})$
be the $(k-1)$-amalgamation of $C$. From the construction of $F_*C'$, it is clear that
$$\newt(F_*C,\eta_{F_*C})=\newt(F_*C',\eta_{F_*C'})\in K^{iso}_1(\bA)~,$$
and by the inductive hypothesis
$$\newt(F_*C',\eta_{F_*C'})=i_*\newt(\bG_*(C'),\eta_{\bG_*(C')}) \in K^{iso}_1(\bA)~.$$
By Proposition \ref{filprop} (iv) a filtered
contraction $\Gamma \colon 1 \simeq 0\colon C \to C$ of $(F_*C,\eta_{F_*C})$ determines a
contraction
$$\bG_*(\Gamma)\colon 1~\simeq~0\colon \bG_*(C) \to \bG_*(C)$$
of the associated $k$-dimensional chain complex
$\bG_*(C)$ in the signed derived category $\SD(\bA)$
with
$$G_r(\Gamma)=(-)^{s}\Gamma_{-1}\colon G_r(C)_s=
C_{r+s,r} \to G_{r+1}(C)_s=C_{r+s+1,r+1}~.$$
In particular, the chain maps
$$\partial_k= d_*\colon G_k(C) \to G_{k-1}(C)~,~
\Gamma_{-1}\colon G_{k-1}(C) \to G_k(C)$$
are related by a chain homotopy
$$\Gamma_0\colon \Gamma_{-1} \partial_k~\simeq~1\colon G_k(C) \to G_k(C)~.$$
Now Proposition \ref{derivedsplit} gives a direct sum system in $\SD(\bA)$
$$\xymatrix@C+5pt{G_k(C)\ \ar@<-.6ex>@{<-}[r]_-{\Gamma_{-1}}
\ar@<.6ex>[r]^-{\partial_k} &
\ G_{k-1}(C)\
\ar@<-.6ex>@{<-}_-{\Delta}[r]
\ar@<.6ex>[r]^-{g} &\  C(\partial_k)}\ . $$
This defines a $k$-splitting of $\bG_*(C)$, and by Lemma \ref{fiveten}
the associated amalgamated complex $\bG_*(C')$ is the corresponding abelian folding. By Proposition \ref{elementary}
there is defined an isomorphism of chain complexes in $\SD(\bA)$
$${\spt{i\choose j}}\colon \bG_*(C)~\cong~\bG_*(C') \oplus \bE$$
where $\bE$ is the elementary complex in $\SD(\bA)$ given by
$G_k(C) \xrightarrow{1} G_k(C)$ in degrees $(k,k-1)$. Since
$\eta_{\bG_*(C)} = \eta_{\bG_*(C') \oplus \bE}$,
the formula in Proposition  \ref{elementary}
 for $\tau^{iso}(\spt{i\atop j})$ gives
$$
i_*\newt(\spt{i\atop j})
=~(-)^{k-1}\newt\left ({ g \choose {\Gamma_{-1}}} \colon
G_{k-1}(C) \to
\cC(\partial_k) \oplus G_k(C)\right )
 \in K^{iso}_1(\SD(\bA))~.
$$
By Corollary \ref{tech} this torsion has image
$$i_*\newt(\spt{i\atop j})=0\in K^{iso}_1(\bA)~,$$
so that
$$\begin{array}{ll}
i_*\newt(\bG_*(C))&=~i_*\newt(\bG_*(C')\oplus
\bE)\\[.6ex]
&=~\newt(\bG_*(C')) \in K^{iso}_1(\bA)
\end{array}$$
which gives the inductive step
$$\begin{array}{ll}
\newt(F_*C,\eta_{F_*C})&=~\newt(F_*C',\eta_{F_*C'})\\[.6ex]
&=~i_*\newt(\bG_*(C'),\eta_{\bG_*(C')})\\[.6ex]
&=~i_*\newt(\bG_*(C),\eta_{\bG_*(C)}) \in K^{iso}_1(\bA)~.
\end{array}$$
\end{proof}

\subsection{The Invariance Theorem for filtered chain equivalences}
We have already defined the filtered mapping cone $F_*\FC(f)$  of a filtered chain map (see Definition \ref{fil} (vii)), with associated complex
$\bG_*(\FC(f))$ given by
$$G_r(\FC(f))_s=G_r(D)_s \oplus G_{r-1}(C)_{s}~,$$
 but now we need the signed version.

\begin{definition}  The \emph{filtered signed mapping cone}
$(F_*\FC(f),\eta_{F_*\FC(f)})$ of a filtered chain map
 is the $(k+1)$-filtered signed complex
defined by the filtered mapping cone $F_*\FC(f)$
with sign terms
$$\begin{array}{l}
\eta_{\bG_*(\FC(f))}=\eta_{\bG_*(D)\oplus S\bG_*(C)}\\[.6ex]
\hphantom{\eta_{G_*(\FC(f))}~}=~
-\beta(\bG_*(D),S\bG_*(C))+\epsilon(\bG_*(D)_{odd},\chi(S\bG_*(C))+
\eta_{\bG_*(D)}-\eta_{\bG_*(C)}\\[.6ex]
\hskip50pt \in \text{im}(\epsilon\colon K_0(\SD(\bA)) \otimes K_0(\SD(\bA)) \to K_1^{iso}(\SD(\bA)))~,\\[.6ex]
\eta_{G_r(\FC(f))}=\eta_{G_r(D)\oplus G_{r-1}(C)}\\[.6ex]
\hskip50pt \in \text{im}(\epsilon\colon K_0(\bA) \otimes K_0(\bA) \to K_1^{iso}(\bA))~~
(0 \leqslant r \leqslant k+1)~.
\end{array}$$
\hfill\qed
\end{definition}
With these sign conventions, the rearrangement map $\rho\colon
\FC(f) \to \cC(f)$ is a \emph{simple} isomorphism.
\begin{lemma}\
$\newt\left ((\FC(f), \eta_{\FC(f)}) \xrightarrow{\rho} (\cC(f),\eta_{\cC(f)})\right ) =0 \in K_1^{iso}(\bA)$.
\end{lemma}
\begin{proof}
The proof is by induction on $k$. The result is true for $k=0$,
since in that case $F_*\FC(f)=\cC(f)$, with $\FC(f)_{r,0}=D_r$ and $\FC(f)_{r,1}=C_{r-1}$,
and $\eta_{F_*\FC(f)}=\eta_{\cC(f)}$.

 Assume that $k \geqslant 1$, and that the result is true for $(k-1)$.
Let
 $$f'\colon (F_*C',\eta_{F_*C'}) \to (F_*D',\eta_{F_*D'})$$
be the chain equivalence of the $(k-1)$-amalgamations  induced
by $f$, with
$$\newt((\FC(f'),\eta_{\FC(f')}) \to (\cC(f),\eta_{\cC(f)}))=0
\in K^{iso}_1(\bA)$$
by the inductive hypothesis.
By Proposition \ref{mappingcone} the unfiltered complex
$\cC(f)=\cC(f')$ has two iterated mapping cone descriptions: from the signed filtered complexes
$F_*E=\FC(f)$ and $F_*E'=\FC(f')$ respectively, we have signed complexes
$T_{k,0}(E)$ and $T_{k,0}(E')$, and a commutative diagram
$$ \xymatrix{T_{k-2,0}(E) \ar[r]\ar@{=}[d]& T_{k+1,0}(E) \ar[r]\ar[d]&
S^{k-1}T_{k+1,k-1}(E)\ar[d]\\
 T_{k-2,0}(E') \ar[r]& T_{k+1,0}(E') \ar[r]&
S^{k-1}T_{k+1,k-1}(E')}$$

But $\eta_E = \eta_{E'}$, so the rearrangement map $\rho\colon E\to E'$ has same absolute torsion as the following composite (distribute $S$, then
interchange)
$$ \xymatrix@R-7pt{S^{k-1}\left (G_{k-2}(C) \oplus S\big (G_k(D)\oplus G_{k-1}(C)\oplus SG_k(C)\big )\right ) \ar[d]^{\ 1\, \oplus\, S}\\
   S^{k-1}\left (G_{k-2}(C)\oplus SG_k(D)\oplus S\big (G_{k-1}(C)\oplus SG_k(C)\big )\right )\ar[d]^{\ flip\, \oplus\,  1}\\
  S^{k-1}\left (SG_k(D)\oplus G_{k-2}(C)\oplus S\big (G_{k-1}(C)
  \oplus SG_k(C)\big )\right )}
$$
since $T_{k-2,0}(E) = T_{k-2}(E')$. Therefore
$$
\newt(E\to E')=  (-)^{k-1}
i_*\epsilon (G_k(D), G_k(C)\oplus SG_{k-1}(C)\oplus
G_{k-2}(C))  \in K_1^{iso}(\bA)
 $$
from the formulas in \cite[Lemma 7]{ak1}.
By Proposition \ref{mappingcone},
$$\newt(\FC(f) \to\cC(f')) = \newt(E \to E') + i_*\eta_{\bG_*(E')} -
i_*\eta_{\bG_*(E)}$$
where the sign term from the filtered mapping cones is
$$\begin{array}{ll}
\eta_{\bG_*(E)} -
\eta_{\bG_*(E')}&=-\beta(\bG_*(D), S\bG_*(C)) + \epsilon(\bG_{odd}(D), \chi(S\bG_*(C)))\\[1ex]
&\ +\beta(\bG_*(D'), S\bG_*(C')) - \epsilon(\bG_{odd}(D'), \chi(S\bG_*(C')))
\end{array}
$$
We first notice that $\bG_{odd}(D) - \bG_{odd}(D') = G_k(D)$ and
$$i_*\chi(S\bG_*(C)) = i_*\chi(S\bG_*(C'))$$
 since $G_{k-1}(C') =
G_{k-1}(C)\oplus SG_k(C)$. Therefore
$$\begin{array}{ll}
i_*\epsilon(\bG_{odd}(D), \chi(S\bG_*(C))) - i_*\epsilon(\bG_{odd}(D'), \chi(S\bG_*(C')))\\[1ex]
\qquad\qquad =  i_*\epsilon(G_k(D),\chi(S\bG_*(C)))= -i_*\epsilon(G_k(D), \chi(\bG_*(C)))\ .
\end{array}
$$
For the $\beta$-terms, we observe that
$\beta(\bG_*(D), S\bG_*(C))=\beta(\bG_*(D), S\bG_*(C'))$, and compute
$$
\begin{array}{ll}
i_*\beta(\bG_*(D), S\bG_*(C'))-i_*\beta(\bG_*(D'), S\bG_*(C'))\\[1ex]
\qquad=(-)^k i_*\left [ \epsilon(G_k(D), \sum\limits_{j\geq 1}G_{k-2j-1}(C))
- \epsilon(G_k(D), \sum\limits_{j\geq 2}G_{k-2j}(C))
\right ]\ .
\end{array}
$$
It follows that
$$i_*\eta_{\bG_*(E)} -
i_*\eta_{\bG_*(E')}=
 (-)^{k-1}
i_*\epsilon (G_k(D), G_k(C)\oplus SG_{k-1}(C)\oplus
G_{k-2}(C)) $$
so the sign terms cancel and $\newt(\FC(f) \to\cC(f'))=0$ as
required.
\end{proof}
We proceed now to the statement of the Invariance Theorem.
It is immediate from the definitions that the associated complex of the
filtered mapping cone $(F_*\FC(f),\eta_{F_*\FC(f)})$
of a chain map $f\colon (F_*C,\eta_{F_*C}) \to
(F_*D,\eta_{F_*D})$ of $k$-filtered signed complexes in $\bA$
is the mapping cone of the associated chain map
$\bG_*(f)\colon (\bG_*(C),\eta_{\bG_*(C)}) \to (\bG_*(D),\eta_{\bG_*(D)})$
in the signed derived category $\SD(\bA)$, so we have
$$(\bG_*(\FC(f)),\eta_{\bG_*(\FC(f))})= (\cC(\bG_*(f)),\eta_{\bG_*(D)\oplus S\bG_*(C)})~.$$
By Theorem \ref{equivalence}, the filtered chain map $f$ is a filtered chain equivalence
if and only if $\bG_*(f)$ is a chain equivalence, or $\cC(\bG_*(f))$ is contractible.

\begin{theorem}[{\bf Invariance}] \label{invariance2}
The torsion of a filtered chain equivalence
$f\colon (F_*C,\eta_{F_*C}) \to (F_*D,\eta_{F_*D})$ of $k$-filtered
signed complexes is
$$\newt(f)=i_*\newt(\bG_*(f)) \in K^{iso}_1(\bA)~.$$
\end{theorem}
\begin{proof} By definition
$$\newt(\bG_*(f))=\newt(\cC(\bG_*(f))) \in K^{iso}_1(\SD(\bA))$$
so we have the formulas
$$\begin{array}{ll}
i_*\newt(\bG_*(f))&=~i_*\newt(\cC(\bG_*(f)))\\[.6ex]
&=~i_*\newt(\bG_*(\FC(f)))\\[.6ex]
&=~\newt(\FC(f),\eta_{\FC(f)})\\[.6ex]
&=~\newt(\cC(f),\eta_{\cC(f)}) \\[.6ex]
&=~\newt(f) \in K^{iso}_1(\bA)~.
\end{array}$$
\end{proof}

\subsection{The filtered dual complex}
Let $(\bA,\ast)$ be an additive category with involution
(see \cite{ra18}).
We now define an involution on a sub-category of the additive category
$\SD(\bA)$ (following \cite[\S 3]{ak1}).

\begin{definition}\label{fivefifteen}
  Given a signed chain complex $C$ in $\bA$, the
\emph{$n$-dual signed chain complex} is $(C^{n-*}, \eta_{C^{n-*}})$,
where
$$d_{C^{n-*}} = (-1)^r d^*_C \colon C^{n-r} = (C_{n-r})^* \to
 C^{n-r+1}$$
 and the sign
 $$
\eta_{C^{n-*}} = (-)^n(\eta_C)^* + \beta(C,C) + \alpha_n(C) \in
K^{iso}_1(\bA)
$$
 where
$$\alpha_n(C) = \sum_{r\equiv n+2,\, n+3\,  (\text{mod} 4)}
\epsilon(C^r, C^r)
 \in K^{iso}_1(\bA)
\ .$$
 For any integer $n\geq 0$, let $\SPD_n(\bA)$ denote the full
 sub-category of $\SD(\bA)$ whose objects are signed chain
 complexes
 $(C,\eta_C)$ in $\bA$ with $\dim C =n$, such that
 $C^{n-*}$ is chain equivalent to $C$ and $\chi(C) = 0$ for $n$ odd.
\hfill\qed
\end{definition}

 \begin{remark}
 The condition that $\chi(C) = 0$ when $n$ is odd is necessary to ensure that $C\mapsto C^{n-*}$ is an additive functor on $\SPD_n(\bA)$.
 This condition and the requirement that $C^{n-*}$ be chain equivalent to $C$ are both satisfied
 (by duality) for symmetric Poincar\'e
 $n$-complexes of $R$-modules, and these are the main objects of
 interest in the rest of the paper.
 \end{remark}
 \begin{example}
 Let $\bA = \bA(R)$ denote the additive category
of finitely generated based free modules over a ring $R$ with involution. Then $\bA(R)$
has an involution $\ast\colon \bA(R) \to \bA(R)$ (see \cite{ra18}),
and we get an involution on $\SPD_n(R):=\SPD_n(\bA(R))$ given by
\[*\colon C\mapsto C^{n-*}\]
\[*\colon (f\colon C \to D) \mapsto (f^{n-*}\colon D^{n-*}\to C^{n-*})\]
(see \cite{ak1})
\hfill\qed
\end{example}
We will now define the \emph{filtered dual}  $\fildual_*C$ of a
$k$-filtered
complex $F_*C$ (under some filtration assumptions to
ensure that each $G_r(C)$ is an $n$-dimensional chain complex). It
will turn out that the associated complex
$\bG_*(F^{dual}C)$ is the $k$-dual of $\bG_*(C)$. Our main example is
the $k$-filtered chain complex $F_*C(\widetilde E)$ of the total space
of a fibre bundle $F\to  E \to B$, where $\dim B = k$ and $\dim F =n$
(see Section \ref{four}).
\begin{definition}
Let $F_*C$ be a $k$-filtered $(n+k)$-dimensional chain complex in $\bA(R)$.
\begin{enumerate}
\item We say that $F_*C$ is \emph{$n$-admissible} if
\begin{enumerate}
\item[(a)] $C_{r,s} = 0$ unless $0\leqslant s \leqslant k$ and
$0\leqslant r-s \leqslant n$.
\item[(b)] $(G_rC)^{n-*}$ is chain equivalent to $G_rC$.\\
\end{enumerate}

\item We define the \emph{filtered dual} $\fildual_*C$ of $F_*C$ to be the
  $(n+k)$-filtered complex with modules $(\fildual_*C)_{r,s} = C^*_{n+k-r,
  k-s}$ for $0\leqslant s\leqslant k$
and differentials $d^{dual}\colon (\fildual_{*}C)_r \to
(\fildual_{*}C)_{r-1}$ the filtered morphism given by the upper
triangular matrix with
components
$$d^{dual}_j = (-)^{r+s +j(n+r)}d^*_j \colon C^*_{n+k-r, k-s} \to
C^*_{n+k-r+1, k-s+j}~.$$
The associated complex of $\fildual_*C$ will be denoted $\bG_*(\fildual C)$.\\
\item We define the \emph{signed filtered dual} of
  $F_*C$ as follows: Since $F_*C$ is $n$-admissible $\bG_*C$ may be
  considered to lie in $\SPD_n(R)$.  This is a category with involution so
  the dual sign $\eta_{(\bG_*C)^{k-*}} \in K_1^{iso}(\SPD_n(\bA(R))$ is
  defined.  We have an obvious functor $\SPD_n(R) \to \SD(R)$ so we
  may consider $\eta_{(\bG_*C)^{k-*}}$ to lie in $K_1^{iso}(\SD(R))$.
  We define the signed filtered dual of $F_*C$
  to be the signed filtered complex whose underlying filtered
  chain complex is $\fildual_*C$ and with signs:
  \[\eta_{\bG_*(\fildual C)} = \eta_{(\bG_*C)^{k-*}} \in K_1^{iso}(\SD(R))\]
  \[\eta_{G_r(\fildual C)} = \eta_{G_{k-r}(C)^{n-*}} \in K_1^{iso}(R)\]
\end{enumerate}
\hfill\qed
\end{definition}
\begin{remark}
It follows from the condition $(G_rC)^{n-*}$ is chain equivalent to
$G_rC$ that $\chi(F_*C)$ is zero if $n$ is odd.
\end{remark}
\begin{lemma}\label{kdual}
The associated complex $\bG_*(\fildual C)$ of the filtered dual of
  $F_*C$ is the $k$-dual of $\bG_*(C)$ in $\SPD_n(\bA)$.
\end{lemma}
\begin{proof}
We have $G_r(F^{dual}C)_s = (\fildual_*C)_{r+s,r}=
C^*_{n+k-r-s,k-r}$ by definition. The internal differential
$d_0\colon G_r(F^{dual}C)_s \to G_r(F^{dual}C)_{s-1}$ on this term
is given by $(-)^{s}d^*_0$. On the other hand, $G_{k-r}(C)^{n-s} =
C^*_{k-r + n-s, k-r}$, so the dual complex has the same chain
groups as the associated complex. The differential on
$G_{k-r}(C)^{n-s}$ is again $(-)^sd^*_0$, since $\dim G_{k-r}(C) =
n$. The signs on $\bG_*(\fildual C)$ were defined to agree with
those of the $k$-dual of $\bG_*(C)$, and in particular
$G_r(F^{dual}C) = (G_{k-r}C)^{n-*}$ as signed chain complexes, for
$0\leqslant r\leqslant k$.
\end{proof}
\noindent
We define the map
\eqncount
\begin{equation}
\label{filtered_dual_theta}
\theta_{F_*C}\colon C^{n+k-*} \to \fildual_*C
\end{equation}
as the direct sum of the following sign maps on the $(r,s)$-summands of the underlying modules:
\[(-)^{s(n+r+1)}\colon C_{n+k-r,k-s}^* \to C_{n+k-r,k-s}^*\]
\begin{lemma}
 $\theta_{F_*C}$ is a chain equivalence of (unfiltered) complexes.
\end{lemma}
\begin{proof}  We have a commutative diagram
$$\xymatrix@C+10pt{C_{n+k-r,k-s}^*\ \ar[r]^{(-)^{s(n+r+1)}}\ar[d]_{(-)^rd^*_j}
&
\ C_{n+k-r,k-s}^*\ar[d]^{(-)^{r+s +j(n+r)}d^*_j}\\
\qquad C_{n+k-r+1,k-s+j}^*\ \ar[r]^{(-)^{s(n+r)+j(n+r)}}
&\qquad C_{n+k-r+1,k-s+j}^*
}$$
where the vertical maps are the components of the differential on $C^{n+k-*}$ and $\fildual_*C$ respectively.
\end{proof}

\begin{proposition}
\label{filtered_dual}
Let $F_*C$ be an admissible, signed, $k$-filtered chain complex over $\bA(R)$, with $\dim C= n+k$. Then
$\newt(\theta_{F_*C}\colon C^{n+k-*} \to \fildual_*C) =0$.
\end{proposition}
\begin{proof}
The proof will be by induction on
$m(C)$, defined as the
smallest integer $m$ such that $F_{m}C=C$.
Suppose first that $m(C)=0$.  Then $C_r = C_{r,0}$ and $G_r(C) = 0$ unless $r=0$. We still allow non-zero
signs $\eta_{G_r(C)}\in K_1(\bZ)$ if $r>0$.
The map $\theta_{F_*C}\colon  C^{n+k-r} \to C^{n+k-r}$ has absolute torsion
$$\begin{array}{ll}
\newt(\theta)
&=\sum\limits_{r=0}^{n+k} \tau^{iso}((-)^{k(n+r+1)}) + \eta_{\fildual_*C}
- \eta_{C^{n+k-*}}\\[1ex]
&= k\cdot\chi(C_{odd}) + \frac{k}{2}(k+1)\chi(C) +\alpha_{n+k}(C) + \alpha_n(C)  \\[1ex]
&=
\begin{cases}
0\\[-.5ex] \chi(C_{odd})\\[-.5ex] 0 \\[-.5ex] \chi(C_{odd})
\end{cases} +
\begin{cases}
0\\[-.5ex] \chi(C)\\[-.5ex] \chi(C) \\[-.5ex] 0
\end{cases} +
\begin{cases}
0\\[-.5ex] \chi(C_{even})\\[-.5ex] \chi(C) \\[-.5ex] \chi(C_{odd})
\end{cases} \text{\ for\ } k \equiv
\begin{cases}
0\\[-.5ex] 1\\[-.5ex] 2 \\[-.5ex] 3
\end{cases}\Mod 4\\[1ex]
&= 0 \in K_1(\bZ)\ .
\end{array}$$

We now move on to the inductive step. Let $m=M(C) > 0$ for $F_*C$ and let $F_*C'$ be the $m$-amalgamation of $F_*C$, considered as a $k$-filtered complex, with
$$C'_{r,s} = \begin{cases}
C_{r,s} &  s < m-1\\[-.5ex]
C_{r,m-1}\oplus C_{r,m} & s=m-1\\[-.5ex]
0 & s > m-1
\end{cases}
$$
with signs $\eta_{\bG_*(C')} = \eta_{\bG_*(C)}$ and
\[\eta_{G_r(C')} = \begin{cases}\eta_{G_r(C)} &
r\neq m-1\\
\eta_{G_{m-1}(C)\oplus SG_m(C)}
 & r = m-1\end{cases}\]
Note that $M(C') = m-1$, and that $C=C'$ as unfiltered signed chain complexes.
Then we have a commutative diagram:
\[\xymatrix{
C^{n+k-*} \ar[rr]^{\theta_{F_*C}} \ar[rd]_{\theta_{F_*C'}}
  & & \fildual_*C \\
  &\fildual_*C' \ar[ru]_\sigma
}\]
with the map $\sigma\colon \fildual_*C' \to \fildual_*C$ induced by the maps
  on the underlying modules: $\sigma=(-)^{k+1+r}\colon C^*_{r,m} \to C^*_{r,}$, and otherwise $\sigma =1$. By our inductive assumption,
  $\newt(\theta_{F_*C'}) =0$, so it remains to compute
  $\newt(\sigma)$.

 We use the $k$-filtered complexes $F_C$ and $F_*C'$
 to express $C$ and $C'$ as an iterated mapping cones.
 Then
 $$\newt(\sigma) = \tau^{iso}(\sigma) + i_*\eta_{\bG(F^{dual}C)} -
i_*\eta_{\bG(F^{dual}C')} \in K_1(\bZ)=\cy 2
 $$
 where $\tau^{iso}(\sigma)$ is the sum  of the term
 $$\newt((-)^{k+m+n+s+1}\colon G_m(C)^{n-s} \to S(SG_m(C)^{n-s})) = (1+k+m+n)\cdot \chi(G_m(C))$$
 by \cite[Lemma 18]{ak1}, plus
 the torsions of rearrangements of the direct sum
 $$S^{k-m}(G_m(C)^{n-*}\oplus  S(G_{m-1}(C)^{n-*}\oplus S(G_{m-2}(C)^{n-*}
 \oplus \dots )\dots ))$$
 Note that the chain groups $G_m(C)^{n-s} = S(SG_m(C)^{n-s})$, although these chain complexes have different signs on their differentials. Let
 $X = SG_m(C)$, $Y = G_{m-1}(C)$ and
  $Z =S(G_{m-2}(C)^{n-*} \oplus \dots )\dots )$.
  Then we rearrange by
  $$\begin{array}{ll}
  SX^{n-*}\oplus S(Y^{n-*} \oplus Z) &\xrightarrow{a} S(X^{n-*}\oplus Y^{n-*} \oplus Z)\\[1ex]
 \xrightarrow{b} S((X\oplus Y)^{n-*}\oplus Z)
 & \xrightarrow{c} S((Y \oplus X)^{n-*}\oplus Z)
  \end{array}
  $$
  The formulas of \cite[Lemma 7]{ak1} give
  $$\newt(a) + \newt(b) + \newt(c) =
  \epsilon(G_m(C), G_{m-2}(C)\oplus \dots \oplus G_0(C))\ .$$
  Now the sign term $i_*\eta_{\bG(F^{dual}C)} -
i_*\eta_{\bG(F^{dual}C')}$ is the sum of two terms
$$
\alpha_k(\bG_*(C)) -  \alpha_k(\bG_*(C'))
=(1+k+m+n)\cdot \chi(G_m(C))$$
and
$$\beta(\bG_*(C), \bG_*(C)) - \beta(\bG_*(C'), \bG_*(C'))
=
\epsilon(G_m(C), G_{m-2}(C)\oplus \dots \oplus G_0(C))
$$
 so $\newt(\sigma)=0$. Hence
$\newt(\theta_{F_*C}) = \newt(\theta_{F_*C'}) =
0$ as required.
\end{proof}

\begin{example}[Tensor Products]
The following special case is used in deriving the product formula for the absolute torsion of symmetric Poincar\'e spaces (see Proposition \ref{product}).
\begin{lemma}
\label{product_dual}
Let $C$ and $D$ be chain complexes of dimension  over $\bA(R)$ and $(\bA(S))$
respectively, where $R$ and $S$ are rings with involution.
Let $\dim C = k $ and $\dim D = n$. Then
we have a chain equivalence
\[\theta_{C\otimes D}\colon (C\otimes D)^{n+k-*} \to C^{k-*}\otimes D^{n-*}\]
given by
\[\theta_{C\otimes D}=(-)^{(k+n+r)(s+n)}\colon C^{k-r}\otimes D^{n-s} \to
C^{k-r} \otimes D^{n-s}\]
Moreover
\[\newt(\theta_{C \otimes D}) =0
\in K_1(R \otimes S)\]
\end{lemma}
\begin{proof}
The associated complex of $F_*(C^{k-*}\otimes D^{n-*})$ is given by
\[\bG_*(C^{k-*}\otimes D^{n-*}):\ \ldots\to C^{k-r}\otimes
D^{n-*}\xrightarrow{(-)^rd_C^*\otimes 1} C^{k-r+1}\otimes D^{n-*}\]
thus $\bG_* (C^{k-*}\otimes D^{n-*}) = \bG_*(C\otimes D)^{k-*}$ considered as the $k$-dual of a chain complex in
 $\SPD_n(\bA(R\otimes S))$.
Hence $F_*(C^{k-*}\otimes D^{n-*}) = F^{dual}_*(C\otimes D)$; moreover
\[\theta_{C\otimes D} = \theta_{F_*(C\otimes D)}\colon (C\otimes D)^{n+k-*}
\to F_*(C^{k-*}\otimes D^{n-*}) = F^{dual}_*(C\otimes D)\]
The result now follows from Proposition \ref{filtered_dual}.
\end{proof}
\end{example}

\bibliographystyle{amsplain}

\begin{thebibliography}{10}

\bibitem{anderson1}
   {D.~R.~Anderson},
    \emph{The {W}hitehead torsion of the total space of a fiber bundle},
    {Topology}
   \textbf{11} (1972), {179--194}.



\bibitem{atiyah1}
M.~F.~Atiyah, \emph{The signature of fibre-bundles}, Global
Analysis (Papers in
  Honor of K. Kodaira), Univ. Tokyo Press, Tokyo, 1969, pp.~73--84.



\bibitem{casson1}
A.~J.~Casson, \emph{Generalisations and applications of block
bundles}, The
  Hauptvermutung book, $K$-Monogr. Math., vol.~1, Kluwer Acad. Publ.,
  Dordrecht, 1996, pp.~33--67.



\bibitem{chs1}
S.~S.~Chern, F.~Hirzebruch, and J.-P.~Serre, \emph{On the index of
a fibered
  manifold}, Proc. Amer. Math. Soc. \textbf{8} (1957), 587--596.



\bibitem{endo1}
H.~Endo, \emph{A construction of surface bundles over
surfaces with
  non-zero signature}, Osaka J. Math. \textbf{35} (1998), no.~4, 915--930.




\bibitem{hrt1}
I.~Hambleton, A.~A.~Ranicki, and L.~R.~Taylor, \emph{Round
${L}$-theory}, J. Pure
  Appl. Algebra \textbf{47} (1987), no.~2, 131--154.



\bibitem{hirzebruch1}
F.~Hirzebruch, \emph{The signature of ramified coverings}, Global
Analysis
  (Papers in Honor of K.~Kodaira), Univ. Tokyo Press, Tokyo, 1969,
  pp.~253--265.



\bibitem{hnk1}
F.~Hirzebruch, W.~D.~Neumann, and S.~S. Koh, \emph{Differentiable
manifolds and
  quadratic forms}, Marcel Dekker Inc., New York, 1971, Appendix II by W.~Scharlau, Lecture Notes in Pure and Applied Mathematics, Vol. 4.



\bibitem{kt1}
S.~Klaus and P.~Teichner, \emph{private communication},  2003.



\bibitem{kodaira1}
K.~Kodaira, \emph{A certain type of irregular algebraic surfaces},
J. Analyse
  Math. \textbf{19} (1967), 207--215.



\bibitem{ak1}
A.~J.~Korzeniewski, \emph{{A}bsolute {W}hitehead torsion}, e-print
math.AT/0502349.



\bibitem{lueck1}
W.~L{\"u}ck, \emph{The transfer maps induced in the
algebraic {$K\sb
  0$}-and {$K\sb 1$}-groups by a fibration. {I}}, Math. Scand. \textbf{59}
  (1986), 93--121.



\bibitem{lueck2}
\bysame, \emph{The transfer maps induced in the algebraic {$K\sb
0$}- and
  {$K\sb 1$}-groups by a fibration. {II}}, J. Pure Appl. Algebra \textbf{45}
  (1987), 143--169.

\bibitem{lueck-ranicki1}
   W.~{L{\"u}ck and A.~A.~Ranicki},
    \emph{Surgery obstructions of fibre bundles},
  {J. Pure Appl. Algebra}
     \textbf{81} (1992),
     {139--189}.


\bibitem{maumary} S.~Maumary, \emph{Contributions \`a la th\'eorie du
type simple d'homotopie}, Comm.  Math.  Helv.  \textbf{44} (1969), 410--437.

\bibitem{meyer1}
W.~Meyer, \emph{Die {S}ignatur von {F}l\"achenb\"undeln},
Math. Ann.
  \textbf{201} (1973), 239--264.



\bibitem{milnor1}
J.~Milnor, \emph{Whitehead torsion}, Bull. Amer. Math. Soc.
\textbf{72} (1966),
  358--426.

\bibitem{munkholm} H.~J.~Munkholm, \emph{Whitehead torsion for PL
fiber homotopy equivalences}, Algebraic Topology, 1978 Waterloo
1978, Springer Lecture Notes \textbf{741}, 90--101 (1979).




\bibitem{neumann1}
W.~D.~Neumann, \emph{Multiplicativity of signature}, J. Pure
Appl. Algebra
  \textbf{13} (1978), 19--31.

\bibitem{pedersen1}
E.~K.~Pedersen, \emph{Geometrically defined transfers, comparisons}, Math. Z. \textbf{180} (1982), 535--544.

\bibitem{putz}
H.~Putz, \emph{Triangulation of fibre bundles}, Canad. J. Math.
\textbf{19} (1967), 499--513.

\bibitem{ra10}
A.~A.~Ranicki, \emph{The algebraic theory of surgery. {I}.
{F}oundations},
  Proc. London Math. Soc. (3) \textbf{40} (1980), 87--192.



\bibitem{ra11}
\bysame, \emph{The algebraic theory of surgery. {I}{I}.
{A}pplications to
  topology}, Proc. London Math. Soc. (3) \textbf{40} (1980), 193--283.



\bibitem{ra13}
\bysame, \emph{The algebraic theory of torsion. {I}.
{F}oundations}, Algebraic
  and geometric topology (New Brunswick, N.J., 1983), Lecture Notes in
  Mathematics, vol. 1126, Springer, Berlin, 1985, pp.~199--237.



\bibitem{ra16}
\bysame, \emph{The algebraic theory of torsion. {I}{I}.
{P}roducts}, $K$-Theory
  \textbf{1} (1987), no.~2, 115--170.


\bibitem{ra18}
\bysame, \emph{Additive ${L}$-theory}, $K$-Theory \textbf{3}
(1989), no.~2,
  163--195.



\bibitem{seifert-threlfall}
H.~Seifert and W.~Threlfall, \emph{Seifert and
{T}hrelfall:  a
  {T}extbook of {T}opology}, Pure and Applied Mathematics, vol.~89, Academic
  Press Inc., New York, 1980.



\bibitem{wall-book}
C.~T.~C. Wall, \emph{Surgery on {C}ompact {M}anifolds}, Academic
Press, 1970. 2nd edition, A.M.S., 1999.

\bibitem{whitehead} J.\ H.\ C. Whitehead,
\emph{Simple homotopy types}, Amer. J. Math. \textbf{72} (1950), 1--57.


\end{thebibliography}

\providecommand{\bysame}{\leavevmode\hbox
to3em{\hrulefill}\thinspace}

\providecommand{\MR}{\relax\ifhmode\unskip\space\fi MR }

\providecommand{\MRhref}[2]{%
  \href{http\colon //www.ams.org/mathscinet-getitem?mr=#1}{#2}
}

\providecommand{\href}[2]{#2}

\end{document}